\documentclass[EJP]{ejpecp} 
\usepackage{xcolor}
\usepackage{hyperref}
\SHORTTITLE{Conservative random walk}
\TITLE{ Conservative random walk
\thanks{J.E.'s research was supported in part by Simons Foundation Grant  579110.}\ 
\thanks{S.V.'s research was  supported in part by Swedish Research Council grant VR~2014-5157 and Crafoord foundation grant 20190667.} }
\AUTHORS{%
  J\'anos Engl\"ander \footnote{Department of Mathematics, University of Colorado at Boulder,  Boulder, CO-80309, USA,
    \EMAIL{janos.englander@colorado.edu}, \url{https://www.colorado.edu/math/janos-englander}}
  \and %% remove this line and below if single author
  Stanislav Volkov \footnote{Centre for Mathematical Sciences, Lund University, Lund 22100-118, Sweden,
    \EMAIL{stanislav.volkov@matstat.lu.se}, \url{http://www.maths.lth.se/\~s.volkov/}}
  }

\KEYWORDS{coin-turning, random walk, conservative random walk, correlated random walk, persistent random walk, Newtonian random walk, recurrence, transience, scaling limits, time-inhomogeneous Markov-processes, Invariance Principle, cooling dynamics, heating dynamics} 

\AMSSUBJ{60G50, 60F05,60J10} % Edit. Separate items with ;
\SUBMITTED{...} % Edit.
\ACCEPTED{...} % Edit.
\VOLUME{0}

\YEAR{2019}

\PAPERNUM{0}

\DOI{10.1214/YY-TN}

\ABSTRACT{Recently, in \cite{EVW2020},  the ``coin-turning walk'' was introduced on $\Z$. It is a non-Markovian process where the steps form a (possibly) time-inhomogeneous Markov chain. In this article, we follow up the investigation by introducing analogous processes in $\Z^d,\ d\ge 2$: at time $n$ the direction of the process is ``updated'' with probability $p_n$; otherwise the next step repeats the previous one. We study some of the fundamental properties of these walks, such as transience/recurrence and scaling limits.

Our results complement previous ones in the literature about ``correlated'' (or ``Newtonian'') and ``persistent'' random walks.
}

\newcommand{\Cov}{\mathsf{Cov}}
\newcommand{\Corr}{\mathsf{Corr}}

\def \e {{\rm e}}
\def \a   {\alpha}

\def \g   {\gamma}

\def \eps {\varepsilon}
\def \P {\mathbb P}

\newcommand{\E}{\mathbb{E}}
\newcommand{\N}{\mathbb{N}}
\newcommand{\Z}{\mathbb{Z}}
\newcommand{\F}{\mathcal{F}}
\newcommand{\R}{\mathbb{R}}
\newcommand{\Mc}{\mathcal{M}}

\newcommand{\1}{\mathbb{1}}

\begin{document}
%%%%%%%%%%%%%%%%%%%%%%%%%%%%%%%%%%%%%%%%%%%%%%%%%%%%%%%%%%%%%%%%%%%
%%                                                               %%
%% No need for \maketitle.                                       %%
%%                                                               %%
%%%%%%%%%%%%%%%%%%%%%%%%%%%%%%%%%%%%%%%%%%%%%%%%%%%%%%%%%%%%%%%%%%%

%%%%%%%%%%%%%%%%%%%%%%%%%%%%%%%%%%%%%%%%%%%%%%%%%%%%%%%%%%%%%%%%%%%
%%                                                               %%
%% Please replace what follows by the body of your article       %%
%% (up to the bibliography):                                     %%
%%                                                               %%
%%%%%%%%%%%%%%%%%%%%%%%%%%%%%%%%%%%%%%%%%%%%%%%%%%%%%%%%%%%%%%%%%%%

\newpage 
\tableofcontents

\section{Introduction}
We are going to study a non-classical random walk. This kind of process has been studied in one-dimension  in~\cite{EVW2020}, and we now define and study the higher dimensional analogs. To avoid ambiguity, in this paper, by geometric distribution we will mean the probability distribution of the number of Bernoulli trials (and not failures) needed to get one success, i.e., the random variable with support on $\{1,2,\dots\}$. It will be denoted by $\mathsf{Geom}(p).$
Finally, by {\it symmetrized geometric distribution} with parameter $p\in(0,1)$ (or $\mathsf{Sgeom}(p)$) we will mean the distribution of a random variable $X$ such that %$\{p_m\}$ with
\begin{equation}\label{symmetrized.geom}
\P(X=m)=p_m=\frac 12 (1-p)^{|m|-1} p,\quad m=\pm 1,\pm 2,\pm3,\dots
\end{equation}
This latter distribution will play an important role during the investigation of the two dimensional homogeneous case; see Section~\ref{SecHom}.
\subsection{The coin-turning process} We start with reviewing the notion of the {\it coin-turning process}, which has recently been introduced in~\cite{EV2018}; note however, that the following definition is slightly different from that in~\cite{EV2018}. 
Let $p_2,p_3, p_4...$  be a given deterministic sequence of numbers such that $p_n\in[0,1]$ for all $n$; define also $q_n:=1-p_n$. We define the following time-dependent ``coin-turning process" $X_n\in\{0,1\}$, $n\ge 1$, as follows. Let $X_1=1$ (``heads") or $=0$ (``tails") with probability $1/2$. For $n\ge 2$,  set recursively
\[
X_n:=\begin{cases} 
0,&\text{with probability } p_n/2;\\
1,&\text{with probability } p_n/2;\\
X_{n-1},&\text{otherwise},
\end{cases}
\]
that is, a fair coin is flipped with probability $p_n$ and stays unchanged  with probability~$q_n$. Additionally, one can define $p_1=1$ and $X_0\equiv 0$.
The process defined above is the same as the one in~\cite{EV2018} with the sequence $\widehat p_n:=p_n/2$, and so turning occurs with probability at most 1/2.

Consider $\overline X_N:= \frac1N \sum_{n=1}^N X_n$, that is, the  empirical frequency of $1$'s (``heads") up to time $N$, in the sequence of $X_n$'s. We are interested in the asymptotic behavior of this random variable when $N\to\infty$. Since we are interested in limit theorems, it is convenient to consider a centered version of the variable $X_n$;  namely $Y_n:=2X_n-1\in\{-1,+1\}$. We have then

\begin{align}\label{Ydrift}
Y_n:=\begin{cases} 
+1,&\text{with probability } p_n/2;\\
-1,&\text{with probability } p_n/2;\\
Y_{n-1},&\text{otherwise}.
\end{cases}
\end{align}
Note that the sequence $\{Y_n\}$ can be defined equivalently as follows: 
\[
Y_n:= (-1)^{\sum_{i=1}^n W_i},
\]
where $W_1,W_2,W_3,...$ are independent Bernoulli variables with parameters $p_1/2,p_2/2,p_3/2,...$, respectively.
Let  also $\F_n:=\sigma(Y_1,Y_2,...,Y_n),\ n\ge 1$.

For the centered variables $Y_n$, we have $Y_j=Y_i(-1)^{\sum_{i+1}^jW_k},\ j>i$, and so, using $\Corr$ and $\Cov$ for correlation and covariance, respectively, one has
\begin{align}\label{eq:eq_eij}
 \Corr(Y_i,Y_j)&=\Cov(Y_i,Y_j)=\E( Y_iY_j)=\E(-1)^{\sum_{i+1}^jW_k}\\
&\ \ \ \ =\prod_{i+1}^j \E (-1)^{W_{k}}=\prod_{k=i+1}^j (1-p_k)
=:e_{i,j};\nonumber\\
 \E (Y_j\mid Y_i)&=Y_i \E(-1)^{\sum_{i+1}^jW_k}=e_{i,j}Y_i.\label{eq:onemore}
\end{align}
The quantity $e_{i,j}$ plays an important role in the analysis of the coin-turning process.

Finally, throughout the paper, we will use the notation $|x|$ for the standard Euclidean norm in~$\R^d$.

\subsection{The one-dimensional coin-turning walk}
We recall from \cite{EVW2020} the definition of the coin-turning walk in one dimension.
\begin{definition}[Coin-turning walk in $\Z$]\label{def:CTW}
The random walk $S$ on $\mathbb Z$ corresponding to the coin-turning, will be called the \emph {coin-turning walk}. Formally, $S_n:=Y_1+...+Y_n$ for $n\ge 1$; we can additionally define $S_0:=0$, so the first step is to the right  or to the left with equal probabilities. As usual, we can then extend $S$ to a {\em continuous time process} by linear interpolation.
\end{definition}

Even though $Y$ is Markovian, $S$ is not. However, the 2-dimensional process $U$ defined by~$U_n:=(S_{n},S_{n+1})$ is Markovian.

In~\cite{EVW2020} the one dimensional coin-turning walk has been investigated from the point of view of transience/recurrence and scaling limits.

\subsection{Generalizing the ``coin-turning walk'' to higher dimensions: ``conservative random walk''}\label{subs: higher.dim}

We define a random walk corresponding to a given sequence $p_n, n=2,3,\dots$ in $[0,1]$ for $d\ge 2$, similarly to the case $d=1$ in Definition \ref{def:CTW}. Now, instead of turning a ``coin'', we have to roll a ``die'' which has $2d$ sides.

The steps  are defined as follows. Let $Y_n\in \{\pm \e_1,\dots,\pm \e_d\}$ where $\e_i$ are the $2d$ unit vectors in $\R^d$, and let $Y_1$ be chosen uniformly from  these  vectors. Let the vectors $Y_1,Y_2,...$ form an inhomogeneous Markov chain with the transition matrix between times $n$ and $n+1$ given by
$$
(1-p_n){\mathbf I}_d+\frac{p_n}{2d}A_d,\ n\ge 1,
$$
where ${\mathbf I}_d$ is the $d\times d$ identity matrix and 
$$A_d:=\begin{bmatrix}
1 & 1 & \dots & 1  \\
1 & 1 & \dots & 1  \\
\vdots &\vdots &\ddots &\vdots\\
1 & 1 & \dots & 1  \\
\end{bmatrix}$$
is the $d\times d$ matrix of ones. 

Now we define the random walk $S$ on $\Z^d$, starting at $z$. Let $S_n:=z+\sum_{i=1}^n Y_i$ for $n\ge 0$ (with the usual convention that $\sum_{i=1}^0=0$),
and denote by $\P_{z}$ the law of this walk. Sometimes we will simply write $\P$ when $z=\mathbf{0}$.
Equivalently, we can define a sequence of independent Bernoulli random variables $\eta_i$, $i=0,1,\dots$, such that $\P(\eta_i=1)=p_i$, and the increasing sequence of stopping times $\tau_j$, such that $\tau_0=0$ and

$$
\tau_{j+1}=\inf\{k>\tau_j: \ \eta_k=1\},\quad j=0,1,2,\dots
$$
At times $\tau_j$ the walk $S_n$ behaves just like a simple symmetric random walk, while in between those times it keeps going in the direction it was going before.

For the sake of completeness, we will include the time-homogeneous case too, that is the case when $p_n=p$ for $n\ge 2$ where $0<p<1$ (when $p=0$ the walk moves in a straight line, while the case $p=1$ corresponds to the classical simple symmetric random walk; so we do not consider these two degenerate cases).

\begin{figure}[htb]
\begin{center}\includegraphics[scale=0.4]{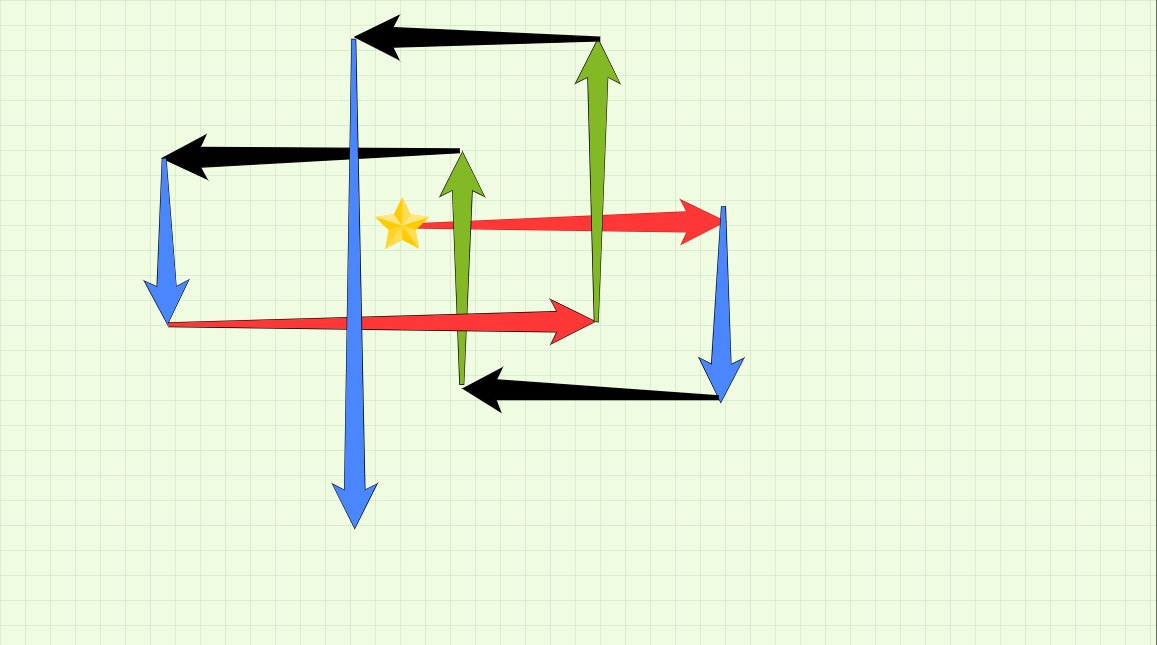}
\caption{A sample path of the walk; the origin is denoted by star. Each arrow keeps track of the ``total'' horizontal/vertical relocation only.}
\label{figEL4}
\end{center}
\end{figure}

Intuitively, the walker is more ``reluctant'' to change direction than an ordinary random walker, motivating the following definition.
\begin{definition}[Conservative random walk]
We dub the process $S$ the {\bf conservative random walk} in $d$ dimensions, corresponding to the sequence $\{p_n\}$.
\end{definition}
\begin{remark} Regarding the sequence of the $p_n$'s we note the following.
\begin{itemize}
\item[(i)] In this paper, we will focus on the case when the $p_n$'s are non-increasing (``cooling dynamics''). Nevertheless,
studying growing $p_n$'s and mixed cases also makes sense. We hope to address this topic in future work.
\item[(ii)] The probability of changing the direction is $p_n\cdot\frac{2d-1}{2d}\le \frac{2d-1}{2d}.$ Our setting thus rules out the kind of heating dynamics (allowed in the setup of ~\cite{EV2018}) when the probability of changing the direction approaches one.$\hfill\diamond$
\end{itemize}
\end{remark}
We now make a fundamental definition.
\begin{definition}[Recurrence/transience]
We call the walk $S$ 
\begin{itemize}
    \item {\bf recurrent} if $\P_{\bf 0}(S_n={\bf 0}\ \text{i.o.})=1$;
    \item {\bf weakly transient}, if it is not recurrent in the above sense;
    \item {\bf strongly transient} if $\P_{\bf 0}(\lim_{n\to\infty}|S_n|=\infty)=1$.
\end{itemize}
\end{definition}
\begin{remark}[Differences compared to traditional categorization]
It is easy to see that strong transience implies weak transience, and that, in fact,
$\P_{z}(\lim_{n\to\infty}|S_n|=\infty)=1$ for each $z\in\Z^d$. 
On the other hand, for recurrence, the probability might depend on the starting point as well as on the ``target.''

Unlike in the case of a simple random walk, it is not {\it a priori} clear whether weak transience necessarily implies strong transience. For example,  the walk might come close to the origin infinitely often, without hitting it, see Figure~\ref{figEL4}. Also, as mentioned above, it is hypothetically possible that the walker visits the origin infinitely often, yet visits some other fixed point only finitely often.

Note that both of these scenarios are possible  only if the probability $p_n$ of updating the direction is not bounded away from zero. Indeed, for a usual random walk, if it hits the origin infinitely often, every time it does, it has a fixed positive probability of, e.g., going to $(1,0)$ on the next step.  Hence, by the usual arguments, it will also hit $(1,0)$ infinitely often. Our random walk, while possibly hitting zero only from a vertical direction, say, at times $\eta_1,\eta_2$, $\eta_3$ etc., might never change its direction at the origin and simply continue going vertically, if $\sum_k p_{\eta_k}<\infty$. Thus, the previous argument will not work.

Finally, although a simple application of Kolmogorov's $0-1$ law  shows (the $\eta_j$ are independent and the  directions chosen at updates too) that $\P_{z}(\lim_{n\to\infty}|S_n|=\infty)\in\{0,1\}$, we cannot rule out the possibility that e.g.
$\P_{\bf 0}(S_n={\bf 0}\ \text{i.o.})\in(0,1).$
$\hfill\diamond$
\end{remark}

\subsection{Some motivation coming from the literature}\label{subs: litera}

Models similar to ours have appeared in the literature under the names ``correlated'' (or ``Newtonian'') random walk and ``persistent random walk.'' (Here ``Newtonian'' refers to the fact that the position of the particle as well as its ``velocity'' affect the next step.) In the statistical physics literature, correlated random walks are often called ``random walks with internal states.''

Scientific phenomena where these models are relevant include polymer growth by sequential addition/deletion of single monomers, and also when a population with births and deaths is considered --- in both cases, growth may predict more growth. Further ones are  flows through a branched structure and the theory of cooperative phenomena in crystals. See \cite{DelaSelvaetal1988,DombFisher1958} for more details on these. 

Following chronological order, we first mention that in \cite{DombFisher1958} random walks with ``restricted reversals'' (i.e. correlated walks) were studied on a class of lattices. Next,
\cite{Gillis1960} treats a two-dimensional problem, where the walker must turn either to the left or to the right, relative to the previous step, with given (constant in time) probabilities.
In \cite{Kac1974} a one dimensional model is considered with a fixed probability of reversal of the last step.

Proceeding to the $1980$s, one dimensional correlated random walks in the homogeneous case are treated in \cite{RenshawHenderson1981}.
Persistent random walks were studied in \cite{SzaszToth1984}, where the ``persistence mechanism'' is given by specifying it at each lattice point and it is done randomly (i.i.d.). In this random environment model, the setup is ``quenched,'' that is, almost sure statements (with respect to the environment) are sought. The main focus of the work was obtaining Central Limit Theorems.
One dimensional correlated random walks are again discussed in \cite{DelaSelvaetal1988} and here even  a related time-inhomogeneous model (``an increasingly more sluggish walk'') is investigated.
The article \cite {Ernst1988} investigates restricted random walks on $d$-dimensional lattices.

Turning to the $1990$s and later, \cite{ChenRenshaw1992,ChenRenshaw1994} revisit the Gillis--Domb-Fisher correlated random walk and generalize it, while \cite{Bohm2000} studies again the one dimensional correlated random walk but this time with two absorbing boundaries. For $d=2$, \cite{Lenci2005} investigates the recurrence of persistent walks. The paper \cite{GruberSchweizer2006} considers (generalized) correlated random walks and studies their diffusive limits. {\em Random flights} are similar processes too. Here a particle in $\mathbb R^d$ changes direction at Poisson times \cite{OrsingherDeGregorio2007a}.

Also, as pointed out to us by Andrew R.\ Wade\footnote{In fact this whole subsection is based on his suggestions.}, correlated random walks may also be discussed by using additive functionals  (this connection is found in \cite{Rogers1985a} and \cite{GeorgiuWade2014}), while some applications in statistical sampling utilize processes ``with momentum'' as well \cite{Bierkensetal2019}.

Finally, recalling that we are mainly interested in the case when $p_n\downarrow 0$, we mention that considering ``cooling'' or (simulated) ``annealing'' is quite standard in the probabilistic/statistical literature. For  papers in this direction which are also quite closely related to our own setup, please see~\cite{Benaimetal,BouguetCloez2018} and the references therein.

\subsection{Outline} The rest of the paper is organized as follows. In Section~\ref{SecHom}, as a warm-up, we prove recurrence when $d=2$ and $p_n=p\in(0,1)$. In Section~\ref{sec_sl2d}, we derive the scaling limit in this case, while in Section~\ref{SecCrit} we do that for the critical case, when the scaling limit is very different 
(a ``zigzag process''). In Section~\ref{Sec2dtrans}, we consider the (recurrent) case when the sequence of the~$p_n$'s is periodic, followed by the proof of transience for the two dimensional walk when $p_n$ has a sufficiently strong decay, as well as that of strong transience for certain multidimensional cases.  In Section~\ref{SecOpen} we formulate some open problems. Finally, the Appendix states and proves some technical lemmas.

\section{Homogeneous case; recurrence on  $\Z^2$}\label{SecHom}
We start with a discussion of the time-homogeneous case. The next result  is perhaps not too surprising.
\begin{theorem}[Recurrence on $\Z^2$; homogeneous case]\label{thm: rechomtwodim} Let $d=2$  and $p_n=p\in (0,1)$, $n\ge 2$. Then  the random walk $S$ is recurrent, that is, $S_n=(0,0)$ infinitely often a.s.
\label{th_recd2}
\end{theorem}
\begin{proof}
(i) Recall that $\tau_1,\tau_2,\dots$ are the consecutive times when the random walk $S$ updates its direction, let $\tau_0=0$ and introduce the {\it embedded walk} $\bar S=(\bar S_k)_{k\ge 0}$ where $\bar S_k:=S_{\tau_{k}}$, $k=0,1,2,\dots$.  This process is a two-dimensional long-range random walk with {\em independent} increments, such that
$$
\bar S_{l+1}-\bar S_{l}=(\kappa_l \xi_l,(1-\kappa_l)\xi_l),
$$
where $\xi_l\sim\mathsf{Sgeom}(p)$ 
as in \eqref{symmetrized.geom},
$\kappa_l$ is Bernoulli($1/2$)
and  $\{\xi_l,\kappa_l\}_l$ is a collection of independent random variables. Equivalently, we can describe $\bar{S}$ as a process with independent increments distributed as
$$
\bar S_{l+1}-\bar S_{l}=\begin{cases}
(|\xi_l|,0), & \text{with probability }1/4;\\
(-|\xi_l|,0), & \text{with probability }1/4;\\
(0,|\xi_l|), & \text{with probability }1/4;\\
(0,-|\xi_l|), & \text{with probability }1/4,
\end{cases}
$$
and the $|\xi_l|$'s are i.i.d.\ \sf{Geom}($p$)  variables. We will show that, in fact, even the embedded process~$\bar S$ is recurrent, and, as a result, so is $S$. To show this, one can directly use Proposition~4.2.4 from~\cite{LL2010} which says that any time-homogeneous random walk on $\Z^2$ with zero drift and finite second moment is recurrent, but for the sake of being self-contained,  we present a short proof based on Lyapunov functions.

Recall that $|(x,y)|=\sqrt{x^2+y^2}$. We use Theorem~2.5.2 from~\cite{MPW2017}, implying that in order to establish recurrence, it is sufficient to find some  function $f:\Z^2\to \R_+$ and~$A>0$ such that
\begin{itemize}
\item $f(z)\to\infty $ as $|z|\to \infty$;
\item if $M_l:=f(\bar S_l)$ then the process $M=\{M_l\}_{l=0,1,2,...}$ satisfies
\begin{align}\label{eqsnos2d}
\E\left[M_{l+1}-M_{l}\mid \bar S_l=(x,y)\right]
\le 0,
\end{align}
whenever $|(x,y)|\ge A$, $\ l=0,1,2,...$
So, informally, $M$ is ``a supermartingale, outside some disc'',
\end{itemize}
We will use the  function $f:\Z^2\to \R_+$ defined as
$$
f(x,y):=\begin{cases}
\ln(x^2+y^2-a)=\ln\left(|(x,y)|^2-a\right), &\text{if } |(x,y)|\ge \sqrt{a+1};\\
0, &\text{otherwise,}
\end{cases}
$$
for some $a>0$ to be chosen later. Denote $r:=|(x,y)|$, $x_1=x+\kappa\xi$ and $y_1=y+(1-\kappa)\xi$, and assume that $r\ge 3$,
where $\xi$ and $\kappa$ have the distribution of $\xi_l$ and $\kappa_l$, respectively. Define the event
$$
{\cal E}=\left \{|\xi|\le r-\sqrt{a+1}\right\}.
$$
By the triangle inequality, on ${\cal E}$ we have $|(x_1,y_1)| \ge \sqrt{a+1}$ and thus $f(x_1,y_1)=\ln(x_1^2+y_1^2-a)$.

Define the random variable $\psi=\psi(x,y)$ for $\|(x,y)\|\ge \sqrt{a+1}$  as
$$
\psi:=\frac{2\xi\zeta +\xi^2}{x^2+y^2-a}
\text{ where } \zeta:=x\kappa+y(1-\kappa)
=\begin{cases}
x, &\text{ with probability }1/2;\\
y, &\text{ with probability }1/2,
\end{cases}
$$
and $\zeta$ is independent of $\xi$. Letting $S_l=(x,y),S_{l+1}=(x_1,y_1)$ a straightforward computation yields that if $r,\|(x_1,y_1)\|\ge \sqrt{a+1}$ then
$$
\Delta_l:=f(\bar S_{l+1})-f(\bar S_{l})=\ln(x_1^2+y_1^2-a)-\ln(r^2-a)=\ln(1+\psi),
$$
Consequently, if $r> \sqrt{a+1}$ then
\begin{align}\label{changed.the.label}
\E\left[\Delta_l\mid \bar S_l=(x,y)\right]&=
\E\left[\Delta_l1_{{\cal E}}\mid \bar S_l=(x,y)\right]+\E\left[\Delta_l1_{{\cal E}^c}\mid \bar S_l=(x,y)\right]\nonumber\\
&=\E\left[\ln(1+\psi) 1_{{\cal E}}\right]+\E\left[\Delta_l1_{{\cal E}^c}\mid \bar S_l=(x,y)\right]\\
&\le\E\left[\ln(1+\psi) 1_{{\cal E}}\right]+\E\left[\ln(x_1^2+y_1^2) 1_{{\cal E}^c}\right]=:(I)+(II).
\nonumber
\end{align}
Another simple computation, using the independence of $\zeta$ and $\xi$, gives
\begin{align}\label{eq:three}
&\E \psi= \frac{\E\xi^2}{r^2-a},
\nonumber\\
&\E \psi^2= \frac{\E\xi^4 +2(x^2+y^2)\E\xi^2 }{(r^2-a)^2},\\
&\E \psi^3= \frac{6(x^2+y^2)\E\xi^4}{(r^2-a)^3}+\mathcal{O}(r^{-6})\E\xi^6,\nonumber
\end{align}
where we used the fact that the odd moments of $\xi$ equal $0$. Also observe that
\begin{align}\label{ximoments}
%\E\xi&=\E\xi^3=\E\xi^5=0,\
\E\xi^2=\frac{2-p}{p^2},\qquad  
\E\xi^4=\frac{(2-p)(p^2+12(1-p))}{p^4},
%\\ \E\xi^6&=\frac{(2-p)(p^4 - 60p^3 + 420p^2 - 720p + 360)}{p^6}
\end{align}
since $\varphi(\lambda)=\E e^{\lambda |\xi|}=\frac{p e^{\lambda}}{1-(1-p)e^{\lambda}}\in(0,\infty)$ for $\lambda<-\log(1-p),$ and $\E \xi^m =\E |\xi|^m=\left.\frac{{\rm d}^m}{{\rm d}\lambda^m}\varphi(\lambda)\right|_{\lambda=0}$ for $0<m\in2\N$.

Now we will use the elementary inequality
\begin{align}\label{logeq}
\ln(1+u)\le u-\frac{u^2}2+\frac{u^3}3\,\qquad \text{for }u>-1,
\end{align}
and observe that (as a brief computation reveals) for fixed $(x,y)$ satisfying $r> \sqrt{a+1}$ we have $\psi>-1$ almost everywhere on the event $\mathcal{E}$. Hence, by~\eqref{eq:three}, \eqref{ximoments} and~\eqref{logeq}, we have that
\begin{align*}
(I)&\le \E\left[ \left(\psi-\frac{\psi^2}2+\frac{\psi^3}3
\right)1_{\cal E}
\right]= \left[\E\psi-\frac{\E(\psi^2)}2+\frac{\E(\psi^3)}3\right]
+(III)
\\ &=-\frac{(2-p)[(2a-3)p^2-36(1-p)] r^2+C_p}{2p^4(r^2-a)^3}+(III)
=-\frac{C+o(1)}{r^4}+(III),
\end{align*}
where 
$$
(III):=-\E\left[ \left(\psi-\frac{\psi^2}2+\frac{\psi^3}3
\right)1_{{\cal E}^c}\right],
$$ 
and $C_p$ is some polynomial of $p$. Note that $C=C(a)>0$  if $2a-3p^2-36(1-p)>0$, that is if we choose  $a=a_p>3/2+\frac{18(1-p)}{p^2}$. We choose such an $a$ and fix it for the rest of the proof.

We also have
$$
|(III)|\le \E\left[ \left(|\psi|+\frac{\psi^2}2+\frac{|\psi|^3}3
\right)1_{{\cal E}^c}\right]\le C_3 e^{-C_2r},
$$
for some $C_2,C_3>0$, since, assuming that $r$, and hence $\max(|x|,|y|)$, is sufficiently large, we get
\begin{align*}
x^2+y^2-a\ge 4\max(|x|,|y|)\ge 2,
\end{align*}
yielding $|\psi|\le \frac12(|\xi|+\xi^2)\le \xi^2$ (recall that $|\xi|\ge 1$) and hence for a positive integer $m$ we have
$\E \left(|\psi|^m 1_{{\cal E}^c}\right)=\E \left(|\psi|^m 1_{|\xi|>r-\sqrt{a+1}}\right)<\E \left(\xi^{2m} 1_{|\xi|>r-\sqrt{a+1}}\right)<C_3' e^{-C_2 r}$ using the properties of the geometric distribution.

Let us also note the  inequalities
\begin{align*}
\ln\left(x_1^2+y_1^2\right)&\le
\ln\left((2x^2+2\kappa^2\xi^2)+(2y^2+2(1-\kappa)^2\xi^2)\right) \\
&=\ln\left(2x^2+2y^2+2\xi^2\right) 
\le \ln (2r^2)+ \ln (2 \xi^2),
\end{align*}
which  use the fact that $\ln(a+b)\le \ln (a)+\ln(b)$ whenever $\min(a,b)\ge 2$.  As a result, the second term in~\eqref{changed.the.label} satisfies that
\begin{align}\label{eqsnos2d2}
(II) &\le \E\left[\left(\ln (2r^2)+ 2\ln |\xi|\right) 1_{{\cal E}^c}\right]
\\ &=  \ln (2r^2)\, \P\left(|\xi|> r-\sqrt{a+1}\right) +2\, \E\left[ \ln\left(2\xi^2\right)
1_{ |\xi|> r-\sqrt{a+1}}\right] \nonumber
 \le C_1 \ln(r) \, e^{-C_2 r}
\nonumber
\end{align}
for some $C_1,C_2>0$, using again % \eqref{easy.bound}.
the properties of the geometric distribution. 

Consequently, from~\eqref{changed.the.label} and~\eqref{eqsnos2d2} we conclude that
\begin{align*}
&\E\left[\Delta_{l}\mid \bar S_l=(x,y)\right] \le -\frac {C+o(1)}{r^4},
\qquad\text{as}\ r\to\infty,
\end{align*}
which is  negative for $r=|(x,y)|$ sufficiently large, and we can apply Theorem~2.5.2 from~\cite{MPW2017}, as alluded to at \eqref{eqsnos2d}, thus completing the proof. 
\end{proof}

\section{Homogeneous case; scaling limit}\label{sec_sl2d}
%In this Section, $d=2$, except in Proposition \ref{prop:LD}. 
We will exploit the following lemma later, but we think that it is also of independent interest.
Let~$L_n$ be a one dimensional coin-turning walk, where $p_n=p\in (0,1)$ for $Y_n$ in~\eqref{Ydrift}.
\begin{proposition}[Tail estimate for the one-dimensional walk]\label{prop:LD}
There exists an $N_0\in\N$ such that if $n\ge N_0$ then for all $a\ge 1$,
\begin{align}\label{fLD}
\P(|L_n|>a \sqrt{n}) \le  f(p,a),\quad\text{where }f(p,a)&:= 2\exp\left\{-p^2\, a/5\right\}.
\end{align}
\end{proposition}
Before presenting the proof of Proposition \ref{prop:LD}, we state and prove a lemma which is a consequence of this proposition.

\begin{lemma}[Upper bound on distance; $d\ge 2$]\label{lemma:estimate}
If $d\ge 2$ then there is an $N_0\in\N$ such that for all $a\ge \sqrt{d}$ and $n\ge N_0$,
\begin{align*}
\P(|S_n|\ge a \sqrt{n}) \le d\, f\left(p,a/\sqrt{d}\right),
\end{align*}
where $f$ is given in~\eqref{fLD}.  
\end{lemma}

\begin{proof}
Let $S^{(j)}$, $j=1,2,\dots,d$  be the $j$-th coordinate of $S$. 
Since after an update, with probability $1/d$ the walk will be moving along the same axis as before, and with probability $1-1/d$ it will start moving in a perpendicular direction, we can write
$$
S^{(j)}_n=L^{(j)}_\kappa
$$
for some $\kappa=\kappa^{(j)}_n\in\{1,2,\dots,n\}$, and
$L^{(j)}$ has the distribution of the one-dimensional walk as in Proposition~\ref{prop:LD}. By Proposition~\ref{prop:LD}, we have
\begin{align*}
\P\left(\left|S^{(j)}_n\right|>a \sqrt{\frac nd}\right)&=
\P\left(\left|L_\kappa^{(j)}\right|>\left(a\sqrt{\frac{n}{d \kappa}}\right) \sqrt{\kappa}\right) \le  f\left(p,a\sqrt{\frac{n}{d \kappa}}\right)
\le  f\left(p,a/\sqrt d\right)
\end{align*}
since $a\mapsto f(p,a)$ is decreasing, 
$n/\kappa\ge 1$, and $a/\sqrt{d}\ge 1$.
This, in turn, implies
\begin{align*}
\P(|S_n|\ge a \sqrt{n})&=\P\left(\sqrt{\sum_{j=1}^d \left(S_n^{(j)}\right)^2}\ge a \sqrt{n}\right)
\le  \P\left(|S_n^{(j)}|\ge a\sqrt {\frac nd}\text{ for some }j\in\{1,2,\dots,d\}\right)
\\ &
\le d\, \P\left(|S_n^{(j)}|\ge a\sqrt {\frac nd}\right)
\le  d \, f\left(p,a /\sqrt{d}\right).
\end{align*} 
\end{proof}

\begin{proof}[Proof of Proposition~\ref{prop:LD}]

Define the strictly increasing integer sequence of stopping times 
$$
0=\tau_0<\tau_1<\tau_2<\dots,
$$
when the walk updates its direction; it keeps going in the same direction between times $\tau_i$ and $\tau_{i+1}$.
Then the $\tau_k-\tau_{k-1}\sim\mathsf{Geom}(p)$, $k=1,2,\dots$ are i.i.d. Moreover,  $\tilde L_k=L_{\tau_k}$ defines the embedded walk, where
$$
\tilde L_k=\xi_1+\dots+\xi_k,
$$
with the i.i.d.\ variables $\xi_i\sim \mathsf{Sgeom}(p)$, while trivially, $|\xi_i|=\tau_i-\tau_{i-1}$.
Let 
\begin{align}\label{eqtau2}
\nu(n):=\min\{j\in\Z_+:\ \tau_{j}\ge n\},
\qquad \text{so that }    \tau_{\nu(n)-1}<n\le \tau_{\nu(n)}
\text{ and }\nu(n)\le n.
\end{align}
Since the walk moves in the same direction between $\tau_{\nu(n)-1}$ and $\tau_{\nu(n)}$, 
\begin{align*}
|L_n| &  \le \max\left(|\tilde L_{\nu(n)-1}|,|\tilde L_{\nu(n)}|\right),
\end{align*}
hence
\begin{align}\label{nunu-1}
\P\left(\left|L_n\right|\ge a\sqrt{n}\right)\le
\P\left(\left|\tilde L_{\nu(n)-1}\right|\ge a\sqrt{n}\right)+
\P\left(\left|\tilde L_{\nu(n)}\right|\ge a\sqrt{n}\right).
\end{align}
Using Markov's inequality we want to bound $\tilde L_m=\xi_1+\dots +\xi_{m}$. Indeed, for any positive integers~$n$, $m\le n$ (see~\eqref{eqtau2}), and any $t\in(0,p)$,
\begin{align}\label{nunu-2}
\P\left(\tilde L_m>a\sqrt{n}\right)
=\P\left(\prod_{i=1}^me^{t\xi_i}>e^{at\sqrt{n}}\right)
\le \frac{\left(\E e^{t\xi_1}\right)^m}{ e^{at\sqrt{n}}}
=e^{-\Lambda(a,t,m,n)},
\end{align}
where
$
\Lambda(a,t,m,n):=a t\sqrt{n}-m\log \left(\E e^{t\xi_1}.
\right).%=a t\sqrt{n}-m\ln\left(\E e^{t\xi_1}\right). 
$
Since
$\E e^{t\xi_1}=
\frac12\E\left(e^{t|\xi_1|}+
e^{-t|\xi_1|}
\right)=\E \cosh({t|\xi_1|})\ge 1,
$
the function $\Lambda$ is monotone decreasing in $m$, hence, for given $p$ and $t$, $\Lambda(a,t,m,n)$ reaches its minimum over $m\le n$ at $m=n$. Note that for $t<-\ln(1-p)$ we have
$$
\E e^{t\xi_1}=\frac p2\left(\frac 1{e^{-t}-(1-p)}+\frac 1{e^{t}-(1-p)}\right).
$$ 
Now it follows via Taylor expansion up to the fourth order of $\Lambda(a,t,n,n)$ with respect to $t$ that
for sufficiently large $n$
\begin{align}\label{eq:Taylor.exp}
\Lambda(a,t(n),n,n)=\frac{p^2(2a-1)}{4-2p}+O(n^{-1})
\ge \frac{p^2\, a}{5} 
\end{align}
where $t(n):=\frac{p^2}{(2-p)\sqrt n}\ (<p<-\ln(1-p))$, since $a\ge 1$ and $p\ge 0$. (Here the term $O(n^{-1})$ is uniform in $a$.)  The result now follows from~\eqref{nunu-1}, \eqref{nunu-2} and \eqref{eq:Taylor.exp}.
\end{proof}

Our next result shows that diffusive scaling  leads to Brownian motion, up to a constant scaling factor.
\begin{theorem}[Scaling limit in the homogeneous case] Let $d\ge 1$ and $p\in(0,1)$. Extend the walk~$S$ to all non-negative times using linear interpolation and for $n\ge 1$, define the rescaled walk ${S}^n$  by 
$${S}^n(t):=\sqrt{\frac{p}{2-p}}\cdot \frac{ S_{nt}}{\sqrt{n}},\ t\ge 0,$$
and finally, let $\mathcal{W}^{(d)}$ denote the $d$-dimensional Wiener measure. Then $\lim_{n\to\infty}\mathsf{Law}({S}^n)=\mathcal{W}^{(d)}$  on $C([0,\infty),\R^d)$. 
\end{theorem}
\begin{remark}
(i) Informally, $ \frac{ S_{n\cdot}}{\sqrt{n}}\approx \sqrt{\frac{2-p}{p}}\cdot B_{\cdot}$, for large $n$, where $B$ is a standard $d$-dimensional Brownian motion. Since $\sqrt{\frac{2-p}{p}}\in (1,\infty)$ for $p\in (0,1)$, the Brownian motion is ``sped up.'' The intuition is that the updates are less frequent compared to a simple random walk, thus there is less cancellation in the steps.

(ii) Note that e.g. for $d=2$, the horizontal and vertical components of the walk are not independent, because, for example, the horizontal component is idle (stays at one location) for the duration of a vertical ``run.''
$\hfill\diamond$\end{remark}
\begin{proof}
%\stas{IS THIS NEEDED? ALREADY IN THE STATEMENT OF THEOREM! As usual for scaling limits of random walks, we automatically  extend all the walks to continuous time processes using linear interpolation, saving us the use of integer parts.}

While using the notation of the previous section, we also take the liberty of using the notation $X_t$ as well as $X(t)$ for a stochastic process $X$, whichever is more convenient at the given instance. We follow the standard route and prove the result by checking the convergence of the finite dimensional distributions (fidis) along with tightness.

\medskip
\noindent \underline{(i) Convergence of fidi's:}
We will argue that the convergence of the fidi's is easy to check for an embedded walk, and the original random walk must have the same limiting fidi's.

To carry out this plan, recall from the proof of Theorem \ref{thm: rechomtwodim}
the  long-range embedded random walk, $\bar S=(\bar S_k)_{k\ge 0}$ where $\bar S_k:=S_{\tau_{k}}$, $k\ge 0$. In this $d$-dimensional setting, its increment vectors are $\bar S_{l+1}-\bar S_{l}= \xi_l\sum_{i=1}^d \1_{\{U_l=i\}}\mathbf{e}_{i},\ l=0,1,...$ where $\xi_l\sim \mathsf{Sgeom}(p)$  and $U_l$ is uniform on $1,2,...,d$ (and the system $\{\xi_l,U_m\}_{l,m\ge 0}$ is independent), while $\{\mathbf{e}_i\}_{1\le i\le d}$ are the  unit basis vectors.
%and that, in fact, even the collection of variables $\{\xi^1_l,\xi^2_l,\dots \xi^d_l\}_{l\ge 0}$ forms an independent system, where $\xi^i_l$ is the $i$th component of $\xi_l,\ 1\le i\le d$.  In particular, unlike in the case of $S$, the  \comj{components of $\bar S$ in different coordinate directions}  are independent.  
Using that $$\mathsf{Var}(\bar S_{l+1}-\bar S_{l})=\E(\bar S_{l+1}-\bar S_{l})^2=\E \xi_l^2 \cdot \E \left|\sum_{i=1}^d \1_{\{U_l=i\}}\mathbf{e}_{i} \right|^2=\E \xi_l^2 \cdot \sum_{i=1}^d \E\1_{\{U_l=i\}} =\frac{2-p}{p^2},$$ it follows that the increment vectors  have  mean value $\mathbf{0}$ and  covariance matrix
$\frac{2-p}{p^2}I_d,$
where $I_d$ is the unit matrix. 
Therefore, denoting $\gamma_p:=\frac p{\sqrt{2-p}}$, we may apply the multidimensional Donsker Invariance Principle (see e.g.\ Theorem 9.3.1 in \cite{Stroock.book}) to the process $\widehat{S}:=\gamma_p \bar{S}$. We obtain that the rescaled walk $\widehat{S}^n$  defined by 
\begin{equation}\label{rescaling}
\widehat{S}^n(t):=\frac{\widehat S_{nt}}{\sqrt{n}}=\gamma_p \frac{\bar S_{nt}}{\sqrt{n}},\ t\ge 0,
\end{equation}  satisfies
\begin{equation}\label{BM.as.limit}
\lim_{n\to\infty}\mathsf{Law}(\widehat{S}^n)=\mathcal{W}^{(d)}\ \text{on}\ C([0,\infty),\R^d).
\end{equation}
In other words,
\begin{equation}\label{rewrite}
\lim_{n\to\infty}\mathsf{Law}\left(t\mapsto \gamma_p \frac{S_{ T(nt)}}{\sqrt{n}}\right)=\mathcal{W}^{(d)}\ \text{on}\ C([0,\infty),\R^d),
\end{equation}
where $T$ is a random time change such that  for  integers  $t=l\ge 0$,  $T(l):=\tau_{l}$ and for $t=l+s,\ s\in (0,1)$, $T(t):=T(l)+(T(l+1)-T(l))s$.

Next,  given that the waiting times for the updates are $\mathsf{Geom}(p)$, the Strong Law of Large Numbers implies that  
\begin{equation}\label{SLLN.T}
\lim_{s\to\infty}T(s)/s\to 1/p,\qquad a.s.
\end{equation} Let $B$ be a standard $d$-dimensional Brownian motion.
We know that for $0\le t_1<t_2<...<t_k$,
\begin{equation}\label{know}
\lim_{n\to\infty}\mathsf{Law}\left(\frac{\gamma_p}{\sqrt{n}}\left(S_{ T(nt_{1})},S_{ T(nt_2)}...,S_{ T(nt_{k})}\right)\right)
=\mathsf{Law}(B_{t_{1}},B_{t_2},...,B_{t_{k}}),
\end{equation}
 and what we want to see next is that this implies that
\begin{equation}\label{almost.there}
\lim_{n\to\infty}\mathsf{Law}\left(
\frac{\gamma_p}{\sqrt{n}}\left( S_{ \frac1pnt_{1}},S_{ \frac1pnt_2},...,
S_{ \frac1pnt_{k}}\right)\right)
=\mathsf{Law}(B_{t_{1}},B_{t_2},...,B_{t_{k}}).
\end{equation}
 It is enough to show (by Slutsky's Theorem) that the difference vector between the vectors on the left-hand sides of \eqref{know} and \eqref{almost.there} converges in probability to $\mathbf{0}$ as $n\to\infty$. (These vectors are such that each of their components are in $\R^d$.)
We will check this component wise. 

Denote by $T^{-1}$ the inverse of the (strictly increasing) map $T$. Fix $\epsilon>0$ and define the random variables $t_{i,n}^*:=\frac{1}{n}T^{-1}\left(\frac 1p nt_{i}\right)$. Then almost surely,
$$\lim_n t^*_{i,n}=\lim_n\frac{1}{n}T^{-1}\left(\frac1 {p}nt_{i}\right)\stackrel{\eqref{SLLN.T}}{=} t_i.$$ 
Fix $1\le i\le k$. The $i$th component  of the difference vector alluded to above, satisfies
$$
\P\left(\frac{\gamma_p}{\sqrt{n}}\left|S_{ \frac1pnt_{i}}- S_{ T(nt_{i})}\right|>\epsilon\right)=\P\left(  \frac{\gamma_p}{\sqrt{n}} \left|S_{T(nt^*_{i,n})}- S_{ T(nt_{i})}\right|    >\epsilon\right),
$$
and we now verify that this converges to zero.  Given that $\lim_n t^*_{i,n}=t_i,$ a.s., it is enough to check that 
$$
\lim_{\delta \to 0}\limsup_{n\to\infty}\P\left(  \frac{\gamma_p}{\sqrt{n}} \left|S_{T(nt^*_{i,n})}- S_{ T(nt_{i})}\right|    >\epsilon\mid E^{i,\delta}_n\right)=0,
$$
where $$E^{i,\delta}_n:=\{ t^*_{i,n}\in (t_i-\delta,t_i+\delta)\}.$$
Since, for any fix $\delta>0$, $\lim_n \P(E^{i,\delta}_n)=1$, there is an $N_\delta\in \N$ such that $\P(E^{i,\delta}_n)\ge 1/2$ for $n>N_{\delta}$,
\begin{align*}
\limsup_{n\to\infty}\,&\P\left(  \frac{\gamma_p}{\sqrt{n}} \left|S_{T(nt^*_{i,n})}- S_{ T(nt_{i})}\right|    >\epsilon\mid E^{i,\delta}_n\right)\le 
\\  &2 \limsup_{n\to\infty}\P\left(\sup_{s\in(t_{i}-\delta,t_{i}+\delta)}\frac{\gamma_p}{\sqrt{n}} \left|S_{T(ns)}- S_{ T(nt_{i})}\right|>\epsilon\right).
\end{align*}
As $\delta\to 0$, the right-hand side tends to zero,
since (as a consequence of the convergence in law to $\mathcal{W}^{(d)}$,)
\begin{align*}
&\lim_n\P\left(\sup_{s\in(t_{i}-\delta,t_{i}+\delta)}\frac{\gamma_p}{\sqrt{n}} \left|S_{T(ns)}- S_{ T(nt_{i})}\right|>\epsilon\right) = 
\mathcal{W}^{(d)}\left(\sup_{s\in(t_{i}-\delta,t_{i}+\delta)}\left|B_{s}- B_{t_{i}}\right|>\epsilon\right).
\end{align*}
We  now have verified \eqref{almost.there}, that is, that 
$$
\lim_{n\to\infty}\mathsf{Law}\left(
t\mapsto
\frac {p}{\sqrt{2-p}}\, 
\frac{S_{ n\tilde t}}{\sqrt{n}}
\right)
=\mathcal{W}^{(d)},$$ where $\tilde t :=\frac{t}{p}$.
Finally, use Brownian scaling: $\tilde t$ can be replaced with $t$, leading to the equivalent limit
$$
\lim_{n\to\infty}\mathsf{Law}\left(t\mapsto\sqrt{\frac{p}{2-p}}\, \frac{ S_{nt}}{\sqrt{n}}\right)=\mathcal{W}^{(d)}.
$$

\noindent\underline{(b) Tightness:}

Exploiting Lemma~\ref{lemma:estimate}, we are going to check Kolmogorov's condition for the fourth moments.  To achieve our goal we fix an $\eps\in(0,p)$ and note that by  Lemma~\ref{lemma:estimate}, and using the fact that $|S_n|\le n$, it follows that
\begin{align*}
\E |S_n|^4=\sum_{i=0}^{n^{4}-1}\P\left(|S_n|^4>i\right)&
\le d^2 n^2
+\sum_{i=d^2 n^2}^{n^{4}-1}\P\left(|S_n|^4>i\right)
\\ & 
=d^2 n^2+\sum_{i=d^2 n^2}^{n^{4}-1} \P\left(|S_n|>a_i\sqrt{n}\right)
\le d^2 n^2+d\sum_{i=d^2 n^2}^{n^{4}-1} f\left(p,a_i/\sqrt{d}\right)
\end{align*}
for all large $n$'s, where $a_i:=\frac{i^{1/4}}{\sqrt{n}}\ge \sqrt d$, as $i\ge d^2 \,n^2$. Since
\begin{align*}\sum_{i=d^2n^2}^{n^{4}-1} f\left(p,a_i/\sqrt{d}\right)=\sum_{i=d^2n^2}^{n^{4}-1} 2e^{-\frac{p^2 a_i}{5\, \sqrt d}}\le 
    \sum_{i=0}^{\infty}2\exp\left\{-\frac{p^2\, i^{1/4}}{5\,\sqrt {d\, n}}\right\},
\end{align*}
by comparing the sum on the right-hand side with the corresponding integral
\begin{align*}
\int_0^\infty   \exp\left\{-\frac{p^2 x^{1/4}}{5\,\sqrt{ d\,n}}\right\}\mathrm{d}x
\stackrel{x=\frac{n^2 d^2 (5u)^4}{p^8} }{=}
\frac{2500\, d^2\, n^2}{p^8} \int_0^\infty u^3  e^{-u}\mathrm{d}u
=\frac{15\,000\, d^2\, n^2}{p^8},
\end{align*} 
we conclude that there exists a $C_p>0$ such that $\E |S_n|^4\le C_p n^2,\ n\ge 1.$
For the rescaled process~$S^{n}$ this yields
$$
\E |S^n(t)|^4=\frac{p^{2}}{{(2-p)}^{2}} \E\left|\frac{ S_{nt}}{\sqrt{n}}\right|^4
\le C_p\frac{n^2t^2}{n^2}=C_pt^2,\ n\ge 1,
$$
(since $\frac{p}{{2-p}}<1$), provided $nt$ is an integer. If $nt$ is not an integer, recall that $S^n$ is defined by linear interpolation and use Jensen's inequality for $y=x^4$ to get the same bound with some $C'_p$ replacing $C_p$.
Finally, by the stationary increments property,
$$
\E |S^n(t)-S^n(s)|^4\le C'_p(t-s)^2,\ 0\le s<t,\ n\ge 1.
$$
 Kolmogorov's condition for tightness is thus satisfied.
\end{proof}

\begin{remark}[A direct bound for $\E|S^4|$, establishing tightness] 
An alternative way of establishing tightness is via computing a bound in a more elementary way for $\E|S_n|^4$.  For simplicity, we will illustrate this in the $d=2$ case. 

Let $L_n:=S^{(x)}_n-S^{(y)}_n$ and $R_n:=S^{(x)}_n+S^{(y)}_n$. Then $L_n=Y_1+\dots+Y_n$ is really a one-dimensional coin-turning walk with parameter $p$ (see~\eqref{Ydrift} and the proof of Lemma \ref{lemma:estimate}.) Next, observe that $|S_n|^2=\left(S^{(x)}_n\right)^2+\left(S^{(y)}_n\right)^2=\frac12(L_n^2+R_n^2)$, yielding that $|S_n|^4\le \frac12\left(L_n^4+R_n^4\right)$. Since $L_n$ and $R_n$ have the same distribution, this implies $\E |S_n|^4\le \E L_n^4$. The latter expectation can be computed directly, albeit that requires a bit of algebra.
%Indeed, we can write $L_n=Y_1+\dots+Y_n$ where $Y_i=(-1)^{W_1+\dots+W_i}$ with the $W_i$ being i.i.d.\ Bernoulli($p/2$). As a result, 
In this time homogeneous case \eqref{eq:eq_eij} and \eqref{eq:onemore} reduce to
\begin{align}\label{YYcond}
e_{i,j}=\Cov(Y_i,Y_j)=\E( Y_iY_j)= q^{j-i};\qquad 
 \E (Y_j\mid Y_i)&=Y_i \E(-1)^{\sum_{i+1}^jW_k}=e_{i,j}Y_i,
\end{align}
for $j\ge i$, where $q:=1-p$.
  We thus have
 \begin{align}\label{ELn4}
 \begin{split}
 &\E L_n^4=\E\left(\sum_1^n Y_i\right)^4\\
 &=\E \left(
  6\sum_{i=1}^{n-1}\sum_{j=i+1}^n Y_i^2 Y_j^2+
 4\sum_{i=1}^{n-1}\sum_{j=i+1}^n Y_i^3 Y_j +
 4\sum_{i=1}^{n-1}\sum_{j=i+1}^n Y_i Y_j^3 +24 \sum_{i=1}^{n-3}\sum_{j=i+1}^{n-2}\sum_{k=j+1}^{n-1}\sum_{l=k+1}^{n} Y_i Y_j Y_k Y_l
 \right. \\ & +\left.
 12\sum_{i=1}^n Y_i^2\left[
 \sum_{k=1}^{i-2}\sum_{l=k+1}^{i-1} Y_k Y_l+
 \sum_{k=i+1}^{n-1}\sum_{l=k+1}^{n} Y_k Y_l+
 \sum_{k=1}^{i-1}\sum_{l=i+1}^{n} Y_k Y_l
 \right]
 +\sum_{i=1}^n Y_i^4 \right)
 \\ &=
 n+3n(n-1)+8 \sum_{i=1}^{n-1}\sum_{j=i+1}^n \E\left[ Y_i Y_j\right]
 +24 \sum_{i=1}^{n-3}\sum_{j=i+1}^{n-2}\sum_{k=j+1}^{n-1}\sum_{l=k+1}^{n}\E\left[ Y_i Y_j Y_k Y_l\right]
 \\ &+
 12\sum_{i=1}^n \left[
 \sum_{k=1}^{i-2}\sum_{l=k+1}^{i-1} \E\left[Y_k Y_l\right]+
 \sum_{k=i+1}^{n-1}\sum_{l=k+1}^{n} \E\left[Y_k Y_l\right]+
 \sum_{k=1}^{i-1}\sum_{l=i+1}^{n} \E\left[Y_k Y_l\right]
 \right].
 \end{split}
 \end{align}
 From~\eqref{YYcond} we obtain that if $i<j<k<l$, then
 %\begin{align*}
 %\E\left[ Y_i Y_j\right]&=
 %\E\left(  Y_i \,\E\left[Y_j| Y_i\right]\right)=
%  \E\left[q^{j-i} Y_i^2\right]=q^{j-i},\\
\begin{align*}
\E\left[ Y_i Y_j Y_k Y_l\right] &
=\E\left(Y_i Y_j Y_k\, \E\left[Y_l|Y_i,Y_j,Y_k\right]\right)
=\E\left(Y_i Y_j Y_k^2 \, q^{l-k} \right)
=q^{l-k}\,\E\left[Y_i Y_j  \right]
=q^{j-i}\times q^{l-k}.
\end{align*}
Substituting this into~\eqref{ELn4} gives
\begin{align*}
\E L_n^4&=n+3n(n-1)
+8\left[\frac{nq}p-\frac{q}{p^2}+O(q^n)\right]
+24\left[\frac{n^2q^2}{2p^2} - \frac{n(5-q)q^2}{(2p^3)} + \frac{3q^2}{p^4}+
+O(n q^n)\right]\\ &
+12\left[
\frac{n^2q}p +\frac{ nq(2q-3)}{p^2} + \frac{2q}{p^2}+O(nq^n)\right]
\\
&= \frac{3n^2(2-p)^2}{p^2}
-\frac{2n(2-p)(p^2+12(1-p))}{p^3}
+
\frac{8(1- p)(3-2p)(3-p)}{p^4}
+O\left(n q^n\right)=O\left(n^2\right),
\end{align*}
hence, Kolmogorov's tightness condition holds.$\hfill\diamond$
\end{remark}

\section{The critical regime}\label{SecCrit}
Next, we turn our attention to the case when $p_n=a/n$ for all large $n$'s where $a>0$. Following \cite{EVW2020}, we call this case the ``critical regime.'' We now need the definition of the ``zigzag process'' in higher dimensions.

\subsection{Preparation:  the zigzag process in higher dimensions}
For simplicity, we start with the two-dimensional case. We describe informally a stochastic process in continuous time, moving in $\R^2$ and starting at the origin. The process is piecewise linear and always moves either horizontally or vertically.

First, notice that if $p^*_n:=\frac{3}{4}p_n$ then the direction is changed with probability $p^*_n$ at time $n$. In our case $p_n=a/n$ and thus $p^*_n=\frac{3a}{4n}=:b/n$ for large $n$'s.

One then takes a realization of the ``scale-free'' Poisson point process, just like it was done in \cite{EVW2020} for $d=1$. This process is defined on $(0,\infty)$ with intensity measure $\frac{b}{x}\, dx$. For the given realization, we construct a trajectory of the zigzag process as follows. Let $t_*$ and $t^*$ in the point process be the left, resp. right neighbors of $t=1$. Toss two independent fair coins and assign one of the labels ``N,W,S,E''  according to the outcome (that is, each has probability $1/4$) to the time interval $(t_*,t^*]$. Going backward in time, label each interval in a way so that the next interval can be labelled in three different ways, each with probability $1/3$, and the label must differ from  that of the previous interval  (if the interval containing $1$ was, say, labeled ``N'', then, going backwards, the next label should be $W$, $S$ or $E$ with equal probabilities, etc.) Do the same for the intervals between the points forward in time. This way, each interval between two consecutive points of the PPP is labeled. All the coin tossings are independent.  These four labels will indicate the direction the process is moving in the time intervals.
 
Let the union of intervals labeled $N$ be $U^{(N)}$ and 
$$
\ell_N(t):=\text{Leb}\left(U^{(N)}\cap[0,t]\right).
$$ 
Define similarly $\ell_S(t),\ell_E(t),\ell_W(t)$.
 
The trajectory of the zig zag process $Z$ will then\footnote{Conditionally on the realization of the PPP and the labeling.} be defined by
$$
Z_t:=(\ell_E(t)-\ell_W(t),\ \ell_N(t)-\ell_S(t)),\qquad t>0.
$$
Clearly, $\lim_{t\to 0} Z_t=(0,0)$ almost surely, even though there are infinitely many points of the PPP in any neighborhood of the origin.
 
When $d>2$, the construction is analogous. The  difference is that in general $p^*_n:=\frac{2d-1}{2d}p_n$ and  $p_n=a/n$ yields $p^*_n=\frac{(2d-1)a}{2dn}=:b/n$ for large $n$'s, and, furthermore, one needs to work with~$2d$ labels. The constraint is that between two consecutive time intervals the process ``must change the label''. 
 
\subsection{Scaling limit for the multidimensional case}
In the critical case, just like in one dimension, proper scaling leads to the zigzag process.
\begin{theorem}[Scaling limit in the critical case]
Let  $d\ge2$ and $p_n=a/n$ for $n\ge n_0.$ 
Extend the walk $S$ to all non-negative times using linear interpolation and for $n\ge 1$, define the rescaled walk ${S}^n$  by 
$$
S^n(t):=\frac{ S_{nt}}{n},\ t\ge 0,
$$
and finally, let $\mathcal{Z}^{(d)}$ denote the law of the $d$-dimensional zigzag process, with parameter $b=\frac{(2d-1)a}{2d}$. Then $\lim_{n\to\infty}\mathsf{Law}({S}^n)=\mathcal{Z}^{(d)}$  on $C([0,\infty),\R^d)$.
\end{theorem} 
\begin{remark}
The result is still valid for $d=1$. Note, however, that the definition of $p_n$  in~\cite{EVW2020} differs by a factor $2$, yielding unit parameter instead of $1/2$.
\end{remark}
\begin{proof} The proof is very similar to that of the one dimensional analog which is part of Theorem~4.11 in~\cite{EVW2020} (see Subsection~6.10 there for the critical case). 
 
 The tightness part  works similarly, namely, just like in \cite{EVW2020}, one simply uses the Lipschitz-1 property of the paths that holds for $S^n$ for each $n\ge 1.$ (This is an advantage compared to the time homogeneous case, and it comes from the fact that the scaling is $n$ and not $n^2$ in this case.)
 
 For the convergence of the finite dimensional distributions, it will be enough to show that weak convergence holds for the processes on $[0,T]$ for any $T>1$.
 
To keep the notation easier, in the rest of the proof we will work with the $d=2$ case, however we note that the general case is completely analogous, by considering $2d$ labels instead of just four.
 
Consider now the space $\mathfrak{S}$ of all double infinite sequences $c_{-2},c_{-1},c_{0},c_1,c_2,...$ where  $c_i\in \{N,W,S,E\}$. When assigning a unique path on $[0,T]$ to a realization of the turning points, the situation is a bit more complicated than for $d=1$. Namely, one has to use a rule, described below in Definition \ref{assigning} with some fixed $s\in\mathfrak{S}$. (In one dimension, there are only two options to assign a path; see Definition 6.10 in \cite{EVW2020}.) Informally, for the segment containing $1$, we assign the label of $c_0$, for the segment to the left and to the right we assign the label of $c_{-1}$  and the label of $c_1$, respectively, etc. 

More precisely, fix $T>1$, denote by ${\Mc}_T$ the set of all locally finite point measures on the interval $(0, T]$, and denote by $N^{(n)}=N^{(n,T)}$ the laws of the point processes induced by the changes of direction of the walk $S^{(n)}$ on the time interval $(0,T]$. 

Let  $s\in\mathfrak{S}$; we now assign a continuous (zigzagged) path to each realization of the point measure.
\begin{definition}[Assigning paths for a given $s\in\mathfrak{S}$]\label{assigning}
Define the map $\Phi_{1}=\Phi_{1,s}: {\Mc}_T \rightarrow C([0,T],\mathbb R^2)$ as follows.
\begin{itemize}
\item First, label the (countably many) atoms on $(0, 1]$ from right to left as $a_1,a_2,...,$ i.e., the closest one on the left to $1$ as $a_1$, the second closest as $a_2$, etc., and note that $1=a_1$ is possible; also label the atoms on $(1, T]$,  from the closest to the furthest as $b_1, b_2$,...;
\item assign label ``$N$"  to all intervals (the union of which is denoted by $S^{(N)}_1$) $[a_i, a_{i+1})$, which are such that in $s$, the corresponding letter is $N$. Here ``corresponding'' means that $c_0$ corresponds to $[a_1,b_1)$, and for $i\ge 1$, $c_i$ corresponds to $[b_i,b_{i+1})$, while $c_{-i}$ corresponds to $[a_{i+1},a_i)$.
\item Do the same for $S, E$ and $W$.
\end{itemize}
The path we obtain will make steps to the North (up) resp. to the West (left), South (down), East (right) on $S^{(N)}_1$, resp. $S^{(W)}_1$, $S^{(S)}_1$, $S^{(E)}_1$.
 
Let  $\mu\in {\Mc}_T$. For  $0<r\leq T$, we define the vertical and horizontal components of the path as
\begin{align*}
\Phi^{\mathsf{vert}}_1(\mu)(r):= \mathsf{Leb}((0,r]\cap S^{(N)}_1)-\mathsf{Leb}((0,r]\cap S^{(S)}_1),\ \mathrm{with}\ \Phi^{\mathsf{vert}}_1(\mu)(0):=0,\\
\Phi^{\mathsf{hori}}_1(\mu)(r):= \mathsf{Leb}((0,r]\cap S^{(E)}_1)-\mathsf{Leb}((0,r]\cap S^{(W)}_1),\ \mathrm{with}\ \Phi^{\mathsf{hori}}_1(\mu)(0):=0,
\end{align*}
and the path itself is $\Phi_1(\mu)(\cdot):=(\Phi^{\mathsf{hori}}_1(\mu)(\cdot),\Phi^{\mathsf{vert}}_1(\mu)(\cdot))$,
where $\mathsf{Leb}$ is the Lebesgue measure on the real line. Then $\Phi^{\mathsf{vert}}_1(\mu)(\cdot)$ is well-defined and continuous on $[0,T]$.  Clearly,
\begin{align} \label{eq:Phi_r}
|\Phi^{\mathsf{vert}}_1(\mu)(r)|,|\Phi^{\mathsf{hori}}_1(\mu)(r)|\leq r,\  0< r\le T.
\end{align}
\end{definition}
Just like in  Proposition 6.13 in \cite{EVW2020}, one can show that for $s\in\mathfrak{S}, T>0$ given,
\begin{itemize}
\item[(a)] $\Phi_{1,s}:{\Mc}_T\to C([0,T],\mathbb R^2)$ is a continuous and uniformly bounded functional, when the former space is equipped with the vague topology, and the latter with  the supremum norm $\|.\|=\|.\|_{[0,T]}$.
\item[(b)]  As $n\to\infty$, $N^{(n)}\rightarrow\sf PPP(b)$ in law (using the vague topology of measures on $(0,T]$), where $\sf PPP(b)$ denotes the law of the Poisson point process on $(0,\infty)$ with intensity measure $b/x\, \mathrm{d}x$.
\end{itemize}

%%% NEW EXPLANATION
Although the proof of (b) can be found in \cite{EVW2020}, in order to be more self contained, we sketch here the main ideas.

Given  $0<a<b<\infty$, $c>0$, set $p_n=\frac{c}{n}\wedge 1,$ and denote the number of updates from step $\lceil an\rceil+1$ to step $\lceil bn\rceil$ by $N^{(n)}((a,b])$. Denoting $\mu_{c;a,b}:=c\ln(b/a)=\int_a^b \frac{c}{x}\, \mathrm{d}x,$ one can show that
\begin{itemize}
\item[(i)] for $k\ge 0,\ 0<a<b$, as $n\to\infty$,
\begin{align}\label{estimation}
\P\left( N^{(n)}((a,b])=k \right)&=\exp(-\mu_{c;a,b}) \frac{\mu_{c;a,b}^k}{k!}+O\left(\frac{1}{n}\right);\\
 {\sf Law}(N^{(n)}((a,b]))&\overset{n\rightarrow\infty}{\longrightarrow} {\sf Poiss}(\mu_{c;a,b});\label{limit_law}
\end{align}
\item [(ii)] given $0< t_1<t_2<...<t_l<\infty$, the random variables $$N^{(n)}( (t_1,t_2]),N^{(n)}((t_2,t_3]),...,N^{(n)}((t_{l-1},t_l])$$ are independent (independent increments), and
\begin{align*}
&{\sf Law}\left( N^{(n)}((t_1,t_2]), N^{(n)}((t_2,t_3]),...,N^{(n)}((t_{l-1},t_l] )\right) \\
&\ \ \ \overset{n\longrightarrow\infty}{\longrightarrow}{\sf PPP(c)} \left( (\mu_{c;t_1,t_2}), (\mu_{c;t_2,t_3})...,(\mu_{c;t_{l-1},t_l}) \right).
\end{align*}
\end{itemize}

One first proves part (i). Once that is done,  since the turns from step $\lceil t_in\rceil +1$ to step $\lceil t_j n\rceil$ and from $\lceil t_l n\rceil+1$ to  $\lceil t_j n\rceil$ are independent for any $0<t_i<t_j\leq t_l<t_r<\infty$, part (ii) will immediately follow.

Regarding part (i), one only needs to prove equation \eqref{estimation}, and then \eqref{limit_law} will easily follow. Checking \eqref{estimation} is done via a straightforward (though a bit tedious) computation. (For $p_n=1/n,   q_n:=1-p_n=(n-1)/n,\ n>>1$, the computation is easier, since  the probability of not updating for a certain time interval then becomes a telescopic product.)

%%% END OF NEW EXPLANATION

After this sketch of the ideas in \cite{EVW2020}, we now return to finish our proof. To accomplish that, let~$P^{\text{uni}}$ be a law on $\mathfrak{S}$ obtained by choosing $C_i$ uniformly in $\{N,W,S,E\}$ and doing it independently for all $i\in\mathcal \Z$, and let   
$$
Q(\cdot):=P^{\text{uni}}(\cdot \mid \forall i\in \Z:\ C_i\neq C_{i+1}).
$$
Furthermore, let $T>1$ and $F:C([0,T],\|\cdot\|)\to\R$ be a bounded continuous functional.  By (a) above,  $F\circ \Phi_{1,s}:(\mathcal{M}_T,\text{vague})\to\R$ is a bounded continuous functional too. Hence, by (b) above, $$\lim_n E((F\circ \Phi_{1,s})(N^{(n)}))=E((F\circ \Phi_{1,s})(\mathsf{PPP}(b))).$$
Finally, bounded convergence  yields that
\begin{equation}\label{exactly.weak.c}
   \lim_n \int_{\mathfrak{S}}E((F\circ \Phi_{1,s})(N^{(n)}))\, Q(\mathrm{d}s)=\int_{\mathfrak{S}}E((F\circ \Phi_{1,s})(\mathsf{PPP}(b)))\, Q(\mathrm{d}s)
\end{equation}
The right-hand side of \eqref{exactly.weak.c} is the expectation of $F$ applied on the zigzag process with parameter~$b$, while the left-hand side is the limit of those terms where the zigzag process is replaced by~$S^{(n)}$ (all processes restricted on $[0,T]$). Since $F$ was an arbitrary bounded continuous functional, \eqref{exactly.weak.c} means that the processes $S^{n}$ converge weakly on the time interval $[0,T]$.
 \end{proof}

%%%%%%%%%%%%%%%%%%%%%%%%%%%%%%%%%%%%%%%%%%%%%%%%%%
%%%%%%%%%%%%%%  NEWER FILE PASTED IN FROM HERE 
%%%%%%%%%%%%%%%%%%%%%%%%%%%%%%%%%%%%%%%%%%%%%%%%%%

\section{Transience/recurrence in higher dimensions\label{Sec2dtrans}}

\begin{figure}[htb]
%\vspace{5mm}
\begin{center}\includegraphics[scale=0.25]{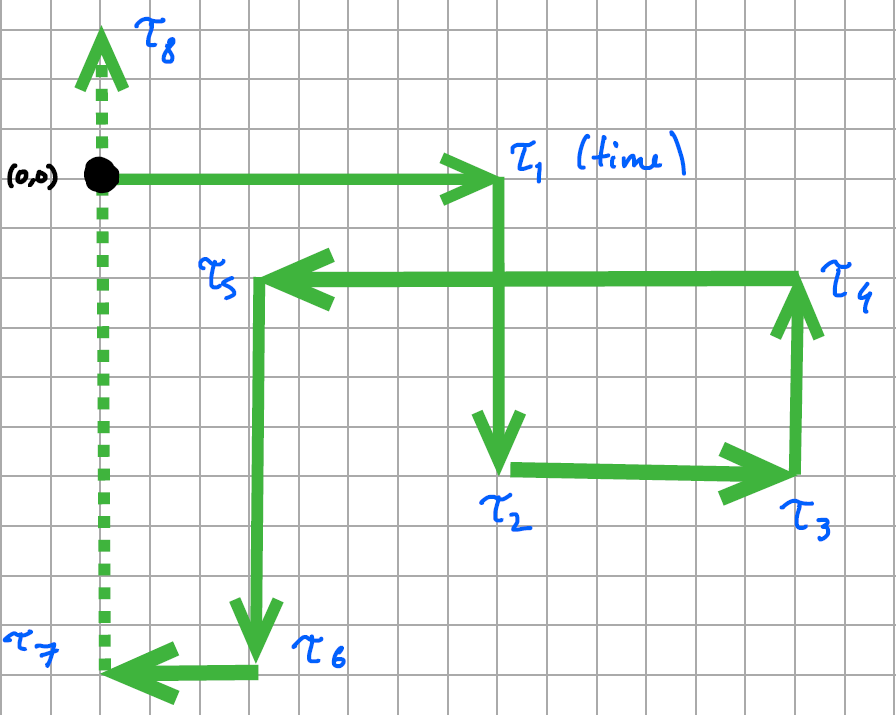}
\caption{A sample path of $S_n$ assuming the walk always turns by $90$ degrees.}
\label{figEL3}
\end{center}
\end{figure}

First, consider the case $d=2$, and define the two-dimensional random walk by $S_n=(X_n,Y_n)\in\Z^2$, $n\ge 1$. The walker keeps going in one of the four directions (up, left, down, or right) at time $n$ with probability $1-p_n$, while with probability~$p_n$ the walker changes direction  in one of the four possible directions, all directions having equal probability.  If $\sum p_n<\infty$, then the walk will make only finitely many turns and then it will trivially drift to infinity. Hence, in the rest of the section we may and will assume $\sum p_n=\infty$, i.e., the walk makes infinitely many turns a.s.

We start with a result that is easy to prove.
\begin{theorem}[Periodic sequence]\label{periodic}
Let the sequence $\{p_n\}$ be periodic, that is, assume that there exists an $r\ge 1$ such that $p_{n+r}=p_n$ for all $n\ge n_0$ with some $n_0\in \N$. Then, for $d=2$ the walk is not strongly transient: $\P(\lim_n|S_n|=\infty)=0$.
\end{theorem}
\begin{proof}
Let $\mathcal{T}\subset \N$ be the set of ``update times.'' 
The proof is based on the following decomposition of the random walk. Let us define a sequence of stopping times: $\tau_0:=n_0$ and 
$$
\tau_{n+1}=\min\{m>n\ :\ r\mid m-n_0,\ m\in\mathcal{T}\}.
$$
Define the random walk $U^1$ on the time interval $\{0,1,2,...,r\}$ as follows
Let 
$$
U^1_i:=S_{n_{0}+i}-S_{n_{0}}
$$ 
and call this a ``block.''
We possibly concatenate further copies $U^2,U^3,...$ in case no update occurs: $U^2$ is added if and only if $n_0+r\not\in \mathcal{T}$, $U^3$ is added if and only if also $n_0+2r\not\in \mathcal{T}$, etc. So the concatenation occurs {\it exactly when} there is no update. Of course, even if there is an update, it might happen (with a chance that is (3/4)th of the update probability) that the direction is unchanged, but we then do not concatenate and the new step is considered to be part of a new~$U^i$ block.

The first step in $U^1$ can be each one of the unit base vectors with equal probabilities, and, by symmetry, this property is inherited for further pieces. The number of pieces (including $U^1$) is geometric with parameter $p_{n_{0}}$. The total length of the walk we obtained this way is exactly~$\tau_1$; this is the first time we have an ``update'' time which is a multiple of $r$. Then repeat the same with the next finite piece of random walk of length $\tau_2-\tau_1$, which is independent of the previous piece, and continue this construction ad infinitum. It is easy to see that the random walk obtained this way is exactly $S$.

Define the embedded walk $S^*$ by $S^*_n:=S_{\tau_{n}}$. The steps of this walk are obtained by concatenating certain blocks, as explained above, and thus they are i.i.d.\ vectors and the length of each one is bounded by the total length of the corresponding piece of the random walk. This latter is a random, geometrically (with parameter $p_{n_{0}}$) distributed multiple of $r$. In particular, the steps of $S^*$ have zero mean and a finite second moment. It follows from Section~8, T1 in~\cite{Spitzer} that $\P(\lim_n|S^*_n|=\infty)=0$ for $d=2$. The same must hold for $S$ too, since $\lim_n|S_n(\omega)|=\infty$ implies $\lim_n|S^*_n(\omega)|=\infty$.
\end{proof}

\begin{theorem}[Strong transience in two dimensions]\label{thm:trans.H.walk}
When $d=2$ and $p_n<n^{-1/2-\eps}$, $n\ge n_0$, for some~$n_0$ and $\eps>0$, the walk is strongly transient, i.e., $| S_n|\to\infty$ a.s.
\end{theorem}
\begin{proof}
First, we prove that $S$ is weakly transient, that is, with probability one it hits $(0,0)$ only finitely often.

This statement is obviously true if $\eps>1/2$, as $\sum p_n<\infty$, so without loss of generality from now on we assume that $\eps\le 1/2$.  Moreover, $p_n<n^{-1/2-\eps'}$ implies $p_n<n^{-1/2-\eps}$ for  $0<\eps<\eps'$; hence, it suffices to prove the theorem only for small positive $\eps$'s. 

As before, let $\tau_n$ be the times when the direction of the walk might change (the $n$th update time); hence, for a fixed $m\ge 1,$
$$
\P(\tau_{n+1}>m+k\mid \tau_n=m,\F_m)=(1-p_{m+1}) (1-p_{m+2})\dots (1-p_{m+k})\stackrel{k\to\infty}{\to} 0
$$
as $\sum p_n=\infty$.
Let us define a subsequence of these stopping times by choosing only those at which the walk switches direction from horizontal to vertical or vice versa. To do so formally, let $\kappa_n\in\{\pm\e_1,\pm\e_2\}$ be the random direction the walk chooses at time $\tau_n$,  set $\eta_0:=0$ and
$$
\eta_{j+1}=\inf\{\N\ni n>\eta_{j}: \ \kappa_{\tau_n}\perp\  \kappa_{\tau_{n-1}}\}.
$$
Also observe that $\eta_{j+1}-\eta_j$ are i.i.d.\  \sf{Geom}($1/2$).

Define the events $A_j$ ($j\ge 1$) as
$$
A_j: =\{\exists a\in \Z\setminus\{{\bf 0}\}:\ X_{\tau_{\eta j}}=(0,a) \text{ or }X_{\tau_{\eta j}}=(a,0)\} ,
$$
that is, $A_j$ is the event that at the  time $\tau_{\eta_j}$ the walker is either on the $x$- or on the $y$-axis. The crucial observation is that in order to hit the origin $(0,0)$ between $\tau_{\eta_j}$ and $\tau_{\eta_{j+1}}$, the event $A_j$ must occur. Indeed, between times $\tau_{\eta_j}$ and $\tau_{\eta_{j+1}}$ the walk moves only along the same line, either $(x,t)$ or $(t,y)$, $t\in\Z$, and unless $x=0$ (resp. $y=0$)  the origin cannot be hit. Furthermore, if the walk is already on the horizontal or vertical axis, then it can hit zero only finitely many times\footnote{in fact, bounded by a Geometric($1/4$) random variable} a.s.\ before leaving this axis. We will show that almost surely, only a finite number of the $A_j$ occur, hence the origin is only visited finitely many times, thus proving non-recurrence.

Fix $j$, denote $\eta_j=:\ell$, and without the loss of generality, suppose that at time $\tau_{\ell}$ the walker starts moving horizontally along the line $(t,y)_{t\in\Z}$ for some $y\ne 0$. Let $m\ge 0$ and
$$
B_{m+1}:=\{\eta_{j+1}=\ell+m+1\}=\{\kappa_{\ell+1},\kappa_{\ell+2},\dots, \kappa_{\ell+m}\in \{-\e_1,\e_1\},\kappa_{\ell+m+1}\in \{-\e_2,\e_2\}\}
$$
and note that $\kappa$'s are i.i.d.\ uniform on $\{\pm e_1,\pm e_2\}$  and are independent of all $\eta$'s. We have for all~$x\in \Z\setminus\{0\}$, that almost surely
\begin{equation}\label{eq-bou}
\begin{split}
&\P\left(A_{j+1}\mid  \F_{\tau_{\ell}},B_{m+1}, S_{\tau_{\ell+m}}=(x,y) \right) 
\\
&\ \ =\frac 12
(1-p_{\tau_{\ell+m}+1})
(1-p_{\tau_{\ell+m}+2})
\dots
(1-p_{\tau_{\ell+m}+|x|-1})
p_{\tau_{\ell+m}+|x|}
<p_{\tau_{\ell+m}+|x|}
\\
&\ \ \le
\frac{1}{(\tau_{\ell+m}+|x|)^{1/2+\eps}}
\le \frac{1}{(\tau_{\ell+m})^{1/2+\eps}}
\le \frac{1}{(\tau_{\eta_j})^{1/2+\eps}}
\end{split}
\end{equation}
(if $x=0$ then the probability on the left-hand side is zero) where the factor $1/2$ comes from the fact that at time $\tau_{\ell+m}$ the walk has an option of going towards or away from zero with equal probabilities; in the final two inequalities we assumed that $j$ is sufficiently large (i.e., $j\ge n_0$ where~$n_0$ is as in the statement) and monotonicity of $\tau_k$'s in $k$.

Since the right-hand side of~\eqref{eq-bou} does not depend on $x$, $y$ and $m$, we conclude that
\begin{align}\label{ehw1}
\P(A_{j+1}\mid  \F_{\tau_{\eta_{j}}}) 
\le  \frac{1}{\tau_{\eta_{j}}^{1/2+\eps}},\ a.s.
\end{align}
for large $j$s.

Our next goal is to show that the $\tau_j$ are ``rare'' in the sense that for $n_j:=(j/8)^{\frac 1{1/2-\eps}}$,
\begin{equation}\label{rare}
\tau_j\ge n_j,\ \text{with finitely many exceptions, a.s.}
\end{equation}
Then~\eqref{rare} implies $\sum_j \P(A_{j+1}\mid \F_{\tau_{\eta_j}}) <\infty$ a.s., and thus (by the conditional Borel-Cantelli lemma) that only a finite number of $A_j$s can occur a.s. Indeed, by~\eqref{ehw1},
$$
\sum_{j=n_0+1}^\infty \P(A_{j+1}\mid \F_{\tau_{\eta_j}}) 
\le\sum_{j=n_0+1}^\infty  \frac{1}{\tau_{j}^{1/2+\eps}},\ a.s.
$$
and, using \eqref{rare}, the sum on the right hand side is a.s. finite, as
$$
\sum_j \frac{1}{n_j^{1/2+\eps}}=
\sum_j \frac{1}{(j/8)^{\frac{1+2\eps}{1-2\eps}}}<\infty.
$$
It remains to verify \eqref{rare}, and as we discussed at the beginning, we may (and will) assume that $\eps\in(0,1/4)$. For this, note that
\begin{align}\label{ehw2}
\P(\tau_j<n)=\P(\xi_1+\dots+\xi_n>j)
\end{align}
where the $\xi_j$s are independent Bernoulli random variables with $\P(\xi_k=1)=p_k$.  Using that $1+x<e^x,\ x>0$ along with Markov's inequality,  we obtain that  for some positive constant~$C_1$,
\begin{align*}
\P(\xi_1+\dots+\xi_n>j)&=\P(2^{\xi_1+\dots+\xi_n}>2^j)
\le 2^{-j} \prod_{i=1}^n \E 2^{\xi_i}
=2^{-j} \prod_{i=1}^n \left(1+p_i\right)
\\ &
<C_1 2^{-j} \prod_{i=1}^n \left(1+\frac 1{i^{1/2+\eps}}\right)
<C_1 2^{-j}\exp\left\{\sum_{i=1}^n i^{-1/2-\eps}\right\}\\ &
<C_1 2^{-j}\exp\left\{\int_{i=0}^n x^{-1/2-\eps}\, dx\right\}=C_12^{-j}\exp\left(\frac{n^{1/2-\eps}}{1/2-\eps}\right)\\ &
\le C_1 \exp\left(4n^{1/2-\eps}-j\ln 2\right),
\end{align*}
which, along with ~\eqref{ehw2} and plugging in $n_j=(j/8)^{\frac 1{1/2-\eps}}$, leads to the estimate
$$
\P(\tau_{j}<n_j)<C_1 e^{(1/2-\ln 2)j}< C_1 e^{-0.193 j}.
$$
Thus \eqref{rare} follows by the Borel-Cantelli Lemma. This completes the proof of non-recurrence.

To upgrade the proof to strong transience is straightforward. Fix any point $(x_0,y_0)\in\Z^2$ and redefine the events $A_j$ as
$$
A_j: = 
\{\exists a\in \Z\setminus\{0\}:\ X_{\tau_{\eta_j}}=(x_0,a) \text{ or }X_{\tau_{\eta_j}}=(a,y_0)\}. 
$$
Following  the weak transience proof verbatim, we obtain that a.s.\ the point $(x_0,y_0)$ is visited  only finitely many times. Hence, the same is true for {\em all} points $(x_0,y_0)$ such that $|(x_0,y_0)|\le r$ for a fixed $r\ge 0$.  We conclude that the walk eventually leaves each given bounded set a.s., which completes the proof of the theorem.
\end{proof}

The following statement is more general than the two-dimensional one in Theorem~\ref{thm:trans.H.walk}, as it works for all $d>1$, however, it requires quite strong regularity conditions.%, which are not needed in Theorem~\ref{thm:trans.H.walk}.

\begin{theorem}[Strong transience for $d\ge 2$]\label{thm: square.summable} Let $d\in\{2,3,4,\dots\}$ and consider the $d$-dimensional coin-turning walk $S$ described in Section~\ref{subs: higher.dim} with $S_n=Y_1+\dots+Y_n$. Assume that for some $\eps>\eps'>0$ and $r> 1$, the sequence of $p_n$'s  satisfies the following conditions:
\begin{align}\label{stascond}
\limsup_{n\to\infty}
\frac{\max\limits_{k\in [0,n^{1-\eps'}]} p_{n-k}}{\min\limits_{k\in [0,n^{1-\eps'}]} p_{n-k}}&<r;
\end{align}
\begin{align}\label{stascond3}
\lim_{n\to\infty} \frac{p_n n^{1-\eps}}{\ln n}=\infty;
\end{align}
\begin{align}\label{stascond2}
\sum_{n=1}^{\infty} \left(\frac{p_n}{n^{1-\eps}}\right)^{d/2}<\infty.&
\end{align}
Then $\sum_{n=1}^{\infty} \P(S_n=w)<\infty$ for any $w\in \Z^d$, and  the walk $S$ is thus strongly transient.
\end{theorem}

\begin{example}[Inverse sub-linear decay]\label{ex:inv.sub.lin} Assume that $\g \in(0,1)$, and $p_n=(c+o(1))/n^{\g}$. Under this assumption \eqref{stascond} is automatically satisfied for any $\eps'>0,$ and assumptions~\eqref{stascond3} and \eqref{stascond2} hold too (pick $\eps<\min\{\gamma,1-\gamma\}$). In this case, therefore, the walk exhibits strong transience. Finally, note that the $\g=1$ assumption (critical case; see~\cite{EV2018}, \cite{EVW2020}) produces a  behavior that is dramatically different  from that in the  $\g \in(0,1)$ case; we believe, nevertheless, that strong transience still holds.
\end{example}

\begin{remark}\label{rem: four.remarks} Concerning the assumptions in Theorem \ref{thm: square.summable}, note that
\begin{itemize} 

\item[(i)] Assumption~\eqref{stascond} implies that $p_n>0$ for all sufficiently large $n$.
\item[(ii)]
Assumptions~\eqref{stascond} and~\eqref{stascond2} for $d=2$ imply that $p_n\to 0$ as $n\to\infty$. Indeed, if along a subsequence, $p_{n_{k}}\ge c>0$ for  $\forall k\ge 1$ then $\sum_{n=1}^{\infty} \frac{p_n}{n^{1-\eps}}\ge \frac{c}{r}\sum_{k=1}^{\infty} \frac{n_{k}^{1-\eps'}}{n_{k}^{1-\eps}}=\infty.$
\item[(iii)]
Assumptions~\eqref{stascond}, \eqref{stascond3}, and \eqref{stascond2} always hold if $d\ge 3$ and $\liminf\limits_{n\to\infty} p_n>0$ (by (ii) this is ruled out when $d=2$), although, then strong transience is anything but surprising.
\item[(iv)] For $d\ge 3$, assumption~\eqref{stascond2} is automatically satisfied when $0<\eps'<\eps<1-2/d.\hfill\diamond$
\end{itemize}
\end{remark}

\begin{proof}
Our goal is to show that for any given lattice point $w\in \Z^d$, one has 
\begin{align}\label{ProfVis great}
    \sum_{n\ge 1}\, \P(S_n=w)<\infty.
\end{align}
The proof will proceed in three steps. First, we introduce a sequence of stopping times and consider the embedded process. Secondly, we show that the total length of the steps of the embedded process during a certain time interval is ``not too large" with high probability. Finally, we estimate the probability of hitting a vertex for the ``remainder'' of the embedded process, using the (integral) inversion formula.

\underline{STEP ONE}: 
Recall that $\eta_j\in\{0,1\}$ was the indicator function of the update  occurring at time~$j$; the~$\eta_j$'s are independent with $\P(\eta_j=1)=1-\P(\eta_j=0)=p_j$. We now consider the walk $S_k$ for times $k=n,n-1,n-2,\dots$ backward. Let $\tau_0:=n$ and let $\tau_k$'s be the decreasing sequence of update times; formally
\begin{align*}
    \tau_{k}&= \max\{j<\tau_{k-1}:\ \eta_j=1\}, \ k=1,2,\dots.
\end{align*}
For definiteness, if for some $k$ we have $\eta_j=0$ for $j=0,1,2,\dots,\tau_{k}-1$, then we set $\tau_{k+1}=\tau_{k+2}=\dots=0$. We will estimate the summands in~\eqref{ProfVis great} as follows. Let $m=m(n):=\lfloor n^{1-\eps}\, p_n\rfloor$ and $V_n:=S_n-S_{\tau_m}$. Clearly, a.s.
\begin{align}\label{eq:clear}
\P(S_n=w)=\P(S_n-S_{\tau_m}=w-S_{\tau_m})=:(*),
\end{align}
and note also that $S_{\tau_m}$ and $S_n-S_{\tau_m}$  are independent. Using the bound
\begin{align}
(*)\le \sup_{z\in\Z^d} \P(S_n-S_{\tau_m}=z)=\sup_{z\in\Z^d}  \P(V_n=z),
\end{align}
it is enough to find numbers $\gamma_n$ such that
\begin{align}\label{eq:enough.to.find}
\sup_{z\in\Z^d}\P(V_n=z)<\gamma_n\ \text{a.s. and}\ \sum_{n\ge 1} \gamma_n<\infty.
\end{align}
For that, the distribution of $V_n$ will be handled by inverting its characteristic function, after which some elementary but tedious computations will be carried out to bound the multiple integrals involved.
In fact, for simplicity (and without loss of rigor), we will assume that $n^{1-\eps}\, p_n$  is an integer.

We will assume that $n$ is so large that the ratio in the $\limsup$ in~\eqref{stascond} does not exceed $r$ and show below that, loosely speaking, 
\begin{itemize}
    \item the probability of the update ``does not change significantly'' over the time segment $[n-n^{1-\eps},n]$;
    \item the variables $l_k=\tau_{k-1}-\tau_k$, $k=1,2,\dots,m$ are ``nearly'' i.i.d. $\mathsf {Geom}(p_n)$;
    \item $\P(l_1+\dots+l_m >n^{1-\eps'}/2)$  is ``very small,''
\end{itemize}
where the $l_k$ are defined via
$$
Z_k:=S_{\tau_{k-1}}-S_{\tau_{k}}\in
\begin{cases}
\{(\pm l_k,0),(0,\pm l_k)\},\ d=2,\\
\{(\pm l_k,0,0),(0,\pm l_k,0),(0,0,\pm l_k)\},\ d=3;\\
\dots
\end{cases}
$$
Note that
$$V_n=Z_1+\dots+Z_m$$ 
($m$ depends on $n$).
To be more precise, first note that given $l_1,l_2,\dots, l_{k-1},l_k$, $Z_k$ has the following distribution: \begin{align*}
Z_k=\begin{cases}
(l_k,0)&\text{with probability }1/4;\\
(-l_k,0)&\text{with probability }1/4;\\
(0,l_k)&\text{with probability }1/4;\\
(0,-l_k)&\text{with probability }1/4,
\end{cases}
\end{align*}
if $d=2$, with  similar formulae for $d\ge 3$.

\medskip
\noindent \underline{STEP TWO}:
Denoting
\begin{align}\label{def:A}
    A:=\{l_1+\dots+l_m =n-\tau_m>  n^{1-\eps'}/2\},
\end{align}
we now  show that 
\begin{align}\label{eq:very.small}
\P(A) =o\left( 
 e^{-n^{\eps-\eps'}}\right),\quad \text{as}\ n\to\infty.
 \end{align}
The event $A$ is the same as having less than $m=n^{1-\eps}p_n$ ``updates'' in the time segment $[n-\frac12 n^{1-\eps'},n]$. The probability of each update at those time points is no less than  $p_n/r$ by~\eqref{stascond}, independently of the others. Hence $\P(A)\le \P(W<m)$ where $W\sim \mathsf {Bin}(N,q)$ with $N=\lfloor n^{1-\eps'}/2\rfloor$ and $q=p_n/r\in(0,1/r]$. Since $(Nq)^i/i!$ is an increasing function in $i$ for $i<Nq$ and $m\le Nq$ as well as $m=o(N)$ as $n\to\infty$, we have
\begin{align*}
\P(W<m)&=\sum_{i=0}^{m-1} 
 {\binom{N}{i}}
q^i (1-q)^{N-i}<
\sum_{i=0}^{m-1} \frac{(N q)^i}{i!} e^{-qN(1+o(1))}
\\ &\le 
m\cdot \frac{(N q)^m}{m!} e^{-qN(1+o(1))}
\le e^{m\log N-qN(1+o(1))}
= e^{-qN(1+o(1))}.
\end{align*}
Since $\ln n=o(m)$ by ~\eqref{stascond3}, we conclude that $\ln n\cdot n^{\eps-\eps'}=o(qN)$, proving \eqref{eq:very.small}. 

\medskip
\noindent\underline{STEP THREE}: Recall that our goal is to find a sequence $(\gamma_n)_{n\ge 1}$ that satisfies \eqref{eq:enough.to.find}. To obtain the distribution of $V_n=S_n-S_{\tau_{m}}$, we invert its characteristic function, and discuss the cases $d=2$ and $d\ge 3$ separately. In the sequel, $\bullet$ will denote the usual dot product in $\R^d$.

\medskip
\noindent\underline{First consider the case $d=2$}. Let $h(t_x,t_y):=\left| \E e^{it \bullet V_n}\right|$,  $t=(t_x,t_y)$. Since $V_n$ has a lattice distribution, the inversion formula is particularly simple (see e.g. Chapter~15, Problem~26 in~\cite{FG1997}, or~\cite{Dynkin}): for any $z=(z_x,z_y)\in \Z^2$  we have 
\begin{equation}\label{eq:star}
\begin{split}
(*):=\P(V_n=z)&=\frac{1}{(2\pi)^2}\int_{0}^{2\pi}\int_{0}^{2\pi} e^{-i(t\bullet z)}\E e^{it \bullet V_n} \, \mathrm{d}t_x\, \mathrm{d}t_y\le \frac{1}{(2\pi)^2}\int_{0}^{2\pi}\int_{0}^{2\pi} h(t_x,t_y) \, \mathrm{d}t_x\, \mathrm{d}t_y
\\ &= \frac{1}{\pi^2}\int_{0}^{\pi}\int_{0}^{\pi} h(t_x,t_y) \, \mathrm{d}t_x\, \mathrm{d}t_y .
\end{split}
\end{equation}
To verify the ultimate equality,\ note that $V_n\in \Z^2$, and its distribution is symmetric in both coordinate variables. Hence $h(u,v)=h(-u,v)=h(u,-v)=h(-u,-v)$ and the integrals of $h$ on the squares $[0,\pi]\times[0,\pi],[0,\pi]\times[\pi,2\pi],[\pi,2\pi]\times[0,\pi]$ and $[\pi,2\pi]\times[\pi,2\pi]$ agree.

Since, given $l_k$, $k=1,2,\dots,m$, $Z_k$ is equally likely to be $(0,\pm l_k)$, $(\pm l_k,0)$ independently of everything else, we have
$$
\E [e^{it \bullet V_n}\mid l_1,\dots,l_m]=\prod_{k=1}^m \E [e^{it \bullet Z_k}\mid l_k]=\prod_{k=1}^m \frac{\cos(t_x l_k)+\cos(t_y l_k)}{2}.
$$ 
Consequently,
\begin{align}\label{eqEitv}
h(t_x,t_y)  &=
\left| \E\left(\E [e^{it \bullet V_n}\mid l_1,\dots,l_m]\right)\right|
\le 
 \E\left(\left|\E [e^{it \bullet V_n}\mid l_1,\dots,l_m]\right|\right)
\le 
 \E\left[\prod_{k=1}^m \phi(t;l_k) \right].%\nonumber
\end{align}
where
$$
\phi(t;l_k)=\phi(t_x,t_y;l_k)=\frac{|\cos(t_x l_k)|+|\cos(t_y l_k)|}{2}.
$$
Since $l_k$ is an integer, 
$$
\phi(\pi-t_x,t_y;l_k)
=\phi(t_x,\pi-t_y;l_k)
=\phi(\pi-t_x,\pi-t_y;l_k)
=\phi(t;l_k).
$$
Consequently, from~\eqref{eq:star},
\begin{align}\label{eq:stas}
 (*)\le \frac{1}{\pi^2}\int_{0}^{\pi}\int_{0}^{\pi} \E\left[\prod_{k=1}^m \phi(t,l_k) \right] \, \mathrm{d}t_x\, \mathrm{d}t_y
=
\frac{4}{\pi^2}\int_{0}^{\pi/2}\int_{0}^{\pi/2} \E\left[\prod_{k=1}^m \phi(t;l_k) \right]\, \mathrm{d}t_x\, \mathrm{d}t_y.
\end{align}
To proceed with the estimation, we now need a lemma. (Recall that the lengths $l_k$ are defined for a given $n$.)
\begin{lemma}\label{lem-cos} Given $(p_k)_{k\ge 1}$ satisfying our assumptions, there exists an $N\in\N$ with the following property.
Let $n\ge N$. Let $2\le k\le m=m(n)$. Then, we have $\omega$-wise that 
\begin{align}\label{eqlemcos}
\E (|\cos(s l_k)|\mid l_1,\dots,l_{k-1})
\le \psi(s) \1_{A_{k-1}^c}+\1_{A_{k-1}},\quad \forall s\in [0,\pi/2],
\end{align}
where 
\begin{align}\label{eq:AAk}
A_{k-1}&=\left\{l_1+...+l_{k-1}>\frac12 n^{1-\eps'}\right\}\stackrel{(m\ge k)}{\subseteq} A,
\end{align}
$$
\psi(s)=\psi_n(s)=\max\left(1-\frac{c_1 s^2}{p_n^2} , 1-c_2\right)
$$
and $c_1>0$, $c_2\in (0,1]$ are constants that depend on $r$ only. 
\end{lemma}
\begin{proof}
Since the inequality is trivial when $A_{k-1}$ occurs,  from now on we assume that we are on~$A_{k-1}^c$. For a fixed $k\in\{2,\dots,m\}$ and $j\ge 1$, let $q_j:=\P(l_k=j\mid  l_1,\dots,l_{k-1})$. Then
\begin{align*}
q_j=p_{n-l_1-\dots-l_{k-1}-j}
 \times
\prod_{i=1}^{j-1}(1-p_{n-l_1-\dots-l_{k-1}-i}),
\end{align*}
where the product equals $1$ by definition for $j=1$. Indeed,  $q_j$ equals the probability of turning first at time $n-l_1-\dots-l_{k-1-j}$ and then keeping the same direction for the consecutive $j-1$ steps.
As long as $l_1+\dots+l_{k-1}+i<n^{1-\eps'}$ by~\eqref{stascond} we have that the above $p$'s lie in $[p_n/r,rp_n]$. Consider two cases (I) $p_n\ge \frac1{2r}$, and (II) $p_n< \frac1{2r}$.

In case (I) apply Lemma~\ref{lemNEW} with $M=1$ and so consider only $q_1$, which (for $n$ large) $\ge \frac1{2r}\times r^{-1}$. Hence, choosing $a=a_{\rm (I)}(r)=\frac1{2r^2}$ in the lemma, yields that the left-hand side of~\eqref{eqlemcos} is smaller than or equal to
\begin{align}
\max\left(1-c_1' a_{\rm(I)} s^2 M^2,1-c_2' a_{\rm(I)} \right)
\stackrel{(M=1)}{=}&\max\left(1-\frac{c_1'  s^2}{2r^2},1-\frac{c_2'}{2r^2} \right)\\
\stackrel{(p_n\ge \frac1{2r})}{\le} &\max\left(1-\frac{c_1'  s^2}{8r^4\, p_n^2},1-\frac{c_2'}{2r^2} \right).
\end{align}
In case (II), note that 
$$
q_j\ge (1-r p_n)^{j-1} \,\frac{p_n}r 
$$
as long as $l_1+\dots+l_{k-1}+j<n^{1-\eps'}$; since we are on $A_{k-1}^c$, this is fulfilled, provided $j<\frac{n^{1-\eps'}}2$. Let $M:=\lfloor 1/(2p_n)\rfloor$ and observe that
\begin{align}\label{eq:M}
\frac{1-r^{-1}}{2p_n}\stackrel{(p_n< \frac1{2r})}\le\frac{1-2p_n}{2p_n}=\frac1{2p_n}-1\le  M\le \frac1{2p_n}
=\frac{n^{1-\eps'}}2\cdot \frac{1}{p_n n^{1-\eps'}}
\stackrel{\eqref{stascond3}}{\le} \frac{n^{1-\eps'}}2,
\end{align}
provided that $n$ is large enough. Thus, for $j=1,2,\dots,M$ we have
$$
q_j\ge (1-rp_n)^{M-1}\frac{p_n}{r}
\ge (1-rp_n)^{\frac{1}{2p_n}-1}\frac{p_n}{r}
\ge (r-1)(1-rp_n)^{\frac{1}{2p_n}-1}(r^2\, M)^{-1}
$$
since $p_n\ge \frac{1-r^{-1}}{M}$. An elementary calculation shows that
\begin{align*}
\inf_{p_n\in(0,(2r)^{-1}]} (1-rp_n)^{\frac{1}{2p_n}-1}
&\ge 
\inf_{p_n\in(0,(2r)^{-1}]} (1-rp_n)^{\frac{1}{2p_n}}
=\inf_{y\in(0,1/2]} (1-y)^{\frac{r}{2y}}
=2^{-r},
\end{align*}
since $(1-y)^{1/y}$ is a decreasing function in $y$. Therefore, $q_j\ge \frac{a_{\rm (II)}}M$ where
$a_{\rm (II)}=\frac{r-1}{2^r\, r^2}$.
Therefore, by Lemma~\ref{lemNEW} with $a=a_{\rm(II)}$ the left-hand side of~\eqref{eqlemcos} does not exceed
\begin{align*}
\max\left(1-c_1' a_{\rm (II)} s^2 M^2,1-c_2' a_{\rm(II)} \right)&\stackrel{\eqref{eq:M}}{\le} 
\max\left(1-\frac{c_1' a_{\rm (II)} s^2 (1-r^{-1})^2}{p_n^2},1-c_2' a_{\rm(II)} \right)\\
&\le 
\max\left(1-\frac{c_1' s^2 (r-1)^3}{2^{r-1} r^4p_n^2},1-c_2' \frac{r-1}{2^{r-1}\, r^2} \right).
\end{align*}
Finally, choosing 
$c_1:=c_1' \min\left(\frac1{8r^4},\frac{(r-1)^3}{2^{r-1}\, r^4}\right)$ and
$c_2:=\min\left(\frac{c_2' }{2r^2},\frac{c_2' (r-1)}{2^r\, r^2},1\right)$ 
concludes the proof.
\end{proof}

Let us now return to the proof of Theorem~\ref{thm: square.summable}. Setting $k:=m$, from Lemma~\ref{lem-cos}, it follows  that when $t_x,t_y\in[0,\pi/2]$, 
\begin{align}\label{eqEitv2}
\begin{split}
&\E\left[\prod_{u=1}^m \phi(t,l_u)\right]\le 
\E \left(\E\left[\prod_{u=1}^m \phi(t,l_u)\mid l_1,\dots,l_{m-1}\right] \right)\\ 
&\le 
\E \left(\prod_{u=1}^{m-1} \phi(t,l_u)\ \left[\frac{\psi(t_x)+\psi(t_y)}2\, \1_{A_{m-1}^c}+\1_{A_{m-1}}\right]\right)\\ 
&\le
\frac{\psi(t_x)+\psi(t_y)}2\,  \E\left[\prod_{u=1}^{m-1} \phi(t,l_u) \right]+\P(A_{m-1})\\ 
&=\dots =
\left(\frac{\psi(t_x)+\psi(t_y)}2\right)^m
+\P(A_{m-1})+\P(A_{m-2})+...+\P(A_{1})\\
&\le 
\left(\frac{\psi(t_x)+\psi(t_y)}2\right)^m
+m\cdot \P(A)
=\left(\frac{\psi(t_x)+\psi(t_y)}2\right)^m
+o\left(n^{1-\eps} e^{-n^{\eps-\eps'}}\right),
\end{split}
\end{align}
by induction on $k=m,m-1,\dots,3,2$ for $\E\left[\prod_{u=1}^{k} \phi(t,l_u) \right]$ and using~\eqref{eq:very.small} and~\eqref{eq:AAk}.

This, along with \eqref{eq:stas}, yields
\begin{align}\label{bound.with.integrals.and.o}
(*)\le \frac 4{\pi^2}\int_0^{\pi/2}\int_0^{\pi/2}  \left(\frac{\psi(t_x)+\psi(t_y)}{2}\right)^m 
\mathrm{d}t_x\, \mathrm{d}t_y+o\left(n^{1-\eps} e^{-n^{\eps-\eps'}}\right).
\end{align}
As far as the second summand is concerned, we can ignore it as it is summable. Focusing on the double integral (in the bound for (*)) only, we will split the area of integration into four regions, depending on whether $t_x$ ($t_y$ resp.) is $>$ or $<\bar t:=p_n\,\sqrt{\frac{c_2}{c_1}}< \frac{\pi}{2}$,
the value of $t$ which makes the candidates for the maximum in $\psi(t)$ equal.  We will also assume that $n$ is large enough (since $p_n\to 0$ when $d=2$ by (ii) of Remark \ref{rem: four.remarks}). 
Then
\begin{align*}
(*)&\le \frac 4{\pi^2}\int_0^{\bar t}\int_0^{\bar t}  \left(\frac{\psi(t_x)+\psi(t_y)}{2}\right)^m 
\mathrm{d}t_x\, \mathrm{d}t_y+
\frac 4{\pi^2}\iint_{t\in \left[0,\frac\pi2\right]^2\setminus [0,{\bar t}]^2} \left(\frac{\psi(t_x)+\psi(t_y)}{2}\right)^m \mathrm{d}t_x\, \mathrm{d}t_y
\\ &=:(I)+(II).
\end{align*}
Turning to scaled polar coordinates: $t_x=\sqrt{2} {p_n\rho \cos \theta}$,
$t_y=\sqrt{2}{p_n\rho \sin \theta}$, we have
\begin{align*}
(I)&= \frac 1{\pi^2}\int_{-\bar t}^{\bar t}\int_{-\bar t}^{\bar t}  \left(1-\frac{c_1 t_x^2+c_1 t_y^2}{2 p_n^2}\right)^m \mathrm{d}t_x\, \mathrm{d}t_y \le 
 \frac {2p_n^2}{\pi^2}\int_0^{2\pi}\int_0^{\sqrt{1/c_1}}  \rho\, \left(1-c_1\rho^2\right)^m 
 \mathrm{d}\rho\, \mathrm{d}\theta
 \\
 &=\frac {4p_n^2}{\pi}\cdot \frac{1}{2c_1(m+1)}
 =\frac{2p_n^2}{\pi c_1(m+1)}\le  \frac{2p_n}{\pi c_1\, n^{1-\eps}}.
\end{align*}
On the other hand,
\begin{align*}
(II)\le   \frac 4{\pi^2}\iint_{t\in \left[0,\frac\pi2\right]^2\setminus [0,{\bar t}]^2} \left[\frac{1+(1-c_2)}{2}\right]^m \mathrm{d}t_x\, \mathrm{d}t_y
\le   \frac 4{\pi^2}\iint_{t\in \left[0,\frac\pi2\right]^2} \left[\frac{2-c_2}{2}\right]^m \mathrm{d}t_x\, \mathrm{d}t_y
= \left(1-\frac{c_2}{2}\right)^m.
\end{align*}
Consequently, almost surely,
\begin{equation*}
\P(S_n-S_{\tau_m}=z)=\P(V_n=z)\le (I)+(II)
=\frac{p_n}{\pi c_1 \, n^{1-\eps}}+\left(1-\frac{c_2}2\right)^{p_n n^{1-\eps}}=:\gamma_n.
\end{equation*}
Since the bound is uniform in $z$, we actually obtained that 
\begin{equation}\label{hurrah}
\sup_{z\in\Z^{2}}\P(V_n=z)\le\gamma_n,
\end{equation}
which is summable in $n$ because of~\eqref{stascond3} and~\eqref{stascond2}. Thus \eqref{eq:enough.to.find} is achieved.

\medskip
\noindent\underline{Consider now the case $d\ge 3$}. 
The proof of the $d=2$ case carries through up to formula~\eqref{eq:star}, but now instead of double integrals, one has to deal with multiple ones. Namely, analogously to~\eqref{bound.with.integrals.and.o}, one now obtains
\begin{align*}
\P(V_n=z)\le \left(\frac{2}{\pi}\right)^d  \idotsint\displaylimits_{[0,\pi/2]^d}  \left(\frac{\psi(t_1)+\psi(t_2)+\dots+\psi(t_d)}{d}\right)^m 
\mathrm{d}t_1\, \mathrm{d}t_2\,\dots \mathrm{d}t_d+o\left(n^{1-\eps} e^{-n^{\eps-\eps'}}\right),
\end{align*}
and we split the area of integration in the same way as in the two-dimensional case, using the same~$\bar t$ as in the two-dimensional case. However, since $d\ge 3$, we can no longer assume that $p_n\to 0$, and it might happen that $\bar t\ge \pi/2$; in this case, a multiple integral analogous to (II) is no longer present in the computation, which makes it easier. On the other hand, when that term is present (i.e. when $\bar t<\pi/2$), its estimate is similar to the two-dimensional case:
\begin{align*}
(II)&= \left(\frac{2}{\pi}\right)^d  \idotsint\displaylimits_{t\in \left[0,\pi/2\right]^d\setminus \left[0,{\bar t}\right]^d} 
\left(\frac{\psi(t_1)+\psi(t_2)+\dots+\psi(t_d)}{d}\right)^m 
\mathrm{d}t_1\, \mathrm{d}t_2\dots \mathrm{d}t_d\\
&\le
\left(\frac{2}{\pi}\right)^d
\idotsint\displaylimits_{t\in \left[0,\pi/2\right]^d} \left(\frac{1+\dots+1+(1-c_2)}{d}\right)^m 
\mathrm{d}t_1\, \mathrm{d}t_2\dots \mathrm{d}t_d
=\left(1-\frac{c_2}{d}\right)^m.
\end{align*}
The estimate of the multiple integral defined analogously to $(I)$ however, will be different:  using scaled $d$-dimensional spherical coordinates 
$t_1=p_n\cdot\rho \cos \theta_1$,
$t_2=p_n\cdot\rho \sin \theta_1 \cos\theta_2$,
$t_3=p_n\cdot\rho \sin \theta_1 \sin\theta_2\cos \theta_3$, $\dots$, 
$t_d=p_n\cdot\rho \sin \theta_1 \sin\theta_2\dots \sin \theta_{d-1}$
and exploiting the facts that the expression in the parenthesis in the integrand below is positive in  the cube $[-\bar t,\bar t]^d$, which lies inside the ball centered at the origin with radius $\sqrt{ d/{c_1}}$ (note that $c_2\le 1$), and the Jacobian $J$ of the transform satisfies $|\mathrm{det}(J)|\le \rho^{d-1}$, one arrives at\footnote{Note: one must replace the limits in the integrals by $\pm\pi/2$ if $\bar t\ge \pi/2$.}
\begin{align*}
(I)&= \frac 1{\pi^d}\idotsint\displaylimits_{[-\bar t,\bar t]^d} \left(1-\frac{c_1 t_1^2+c_1 t_2^2+\dots+c_1 t_d^2}{p_n^2\, d}\right)^m 
\mathrm{d}t_1\, \mathrm{d}t_2\dots \mathrm{d}t_d
\\ & \le 
 \frac {p_n^d}{\pi^d}
 \idotsint\displaylimits_{\rho\in \left[0,\sqrt{\frac{d}{c_1}}\right] \theta_1,\dots,\theta_{d-2}\in[0,\pi],\theta_{d-1}\in[0,2\pi]}
  \rho^{d-1}\,  \left(1-\frac{c_1\rho^2}{d}\right)^m  \mathrm{d}\rho\, 
 \mathrm{d}\theta_1\dots \mathrm{d}\theta_{d-1}
 \\ &
 = \frac {2p_n^d}{\pi}\int_0^{\sqrt{ d/{c_1}}}  \rho^{d-1}\,   \left(1-\frac{c_1\rho^2}{d}\right)^m  \mathrm{d}\rho
  =\frac {p_n^d}{\pi} \left(\frac{d}{c_1}\right)^{d/2} 
 \int_0^1  u^{{\frac{d}2-1}}\,  (1-u)^m  \mathrm{d}u  \\ &
  = {\rm const}\  p_n^d\, \cdot \mathsf{B}\left(\frac{d}{2},m+1\right)\sim {\rm const}\  \Gamma(d/2)p_n^d\, \cdot m^{-d/2}\sim
  {\rm const}\  \Gamma(d/2)\left(\frac{p_{n}}{n^{1-\eps}}\right)^{d/2},
  \end{align*} where $\mathsf{B}$ denotes the Beta-function, and we exploited its well known asymptotics, while we also recalled that $m=\lfloor  n^{1-\eps}\, p_n\rfloor$.
By~\eqref{stascond2}, the last expression is summable in $n$.  As a result, the bound $(I)+(II)$ is uniform in $z$ and is summable in $n$, just as in the case $d=2$. Again, \eqref{eq:enough.to.find} is confirmed.
\end{proof}

\section{Some open problems}\label{SecOpen}
In this section we formulate some open problems.
\begin{problem}[Monotonicity] Fix $d\ge 1$. Is it true that if $p_n'\le p_n$ for $n\ge 1$ and the walk exhibits strong transience for the sequence $\{p_n\}$ then the same holds for the sequence $\{p'_n\}$? (Compare with the last sentence in Example~\ref{ex:inv.sub.lin}.)
\end{problem}

\begin{problem}[Critical case] Fix $d\ge 2$. In the critical case ($p_n=\mathrm{const}/n$ for large $n$) we conjecture that the walk exhibits strong transience, which, of course, would readily follow form monotonicity (cf. the last sentence in Example \ref{ex:inv.sub.lin}) and that this property is inherited to its scaling limit, the zigzag process too.
\end{problem}
\begin{problem}[Transient dimensions] Is it true that the walk always exhibits strong transience whenever $d\ge 3$? Clearly, monotonicity would imply this, since when $p_n=p=1,$ $S$ is a simple symmetric random walk.

In the periodic case, for example, using the notation of the proof of Theorem~\ref{periodic}, when~$d\ge 3$, it is known (see, e.g.~\cite{Spitzer}) that  $\P(\lim_n|S^*_n|=\infty)=1.$ We know that for all $\tau_n<m<\tau_{n+1}$  $(1/r)|S_m-S^*_n|$ is uniformly bounded by a geometrically  distributed variable, with parameter~$p_{n_{0}}$. However, it is not clear whether this is enough to control $S$ via $S^*$.

\end{problem}
\begin{problem}[Slow decay]
What happens  when $d=2$ and the $p_n$'s decay slowly as $n\to\infty$? We already know that if $p_n=\text{const}\in (0,1)$ then the walk is recurrent. An interesting question is whether the answer changes if, for example, $p_n\sim 1/\log n$.
\end{problem}

\section{Appendix}
Here we state and prove two lemmas that were needed in the proofs. Some versions of these statements are presumably known, but since we could not find a proper reference, we present their proofs here.
\begin{lemma}\label{lem:arithm}
Let $s\in (0,1/2]$, $s_0\in \R$, and $M\ge 2$ be an integer such that $Ms\ge 1$. Moreover, let $a_k=ks+s_0\  
(\text{mod}\ 1)$, for $k=1,2,\dots,M$. Then\footnote{The constant $2/15$ is definitely sub-optimal. The true constant is probably $1/3$, with equality achieved for $M=3,6,9,\dots$.}
\begin{align*}
N_{s,M}&:=|\{k:\ a_k\in [0,1/2)\}|\ge \frac{2}{15}M.
\end{align*}
\end{lemma}
\begin{proof}
First, assume that $s\ge 1/4$. Since the  length of the interval $[1/2,1)$  is $1/2$, and $s\le 1/2$, for any triple $\{a_i,a_{i+1},a_{i+2}\}$, $i\ge 1$, at least one element must not lie in $[1/2,0)$. Therefore, $N_{s,M}\ge \lfloor M/3\rfloor\ge M/5$.

Next, assume that $s<1/4$, and define
$$
\tau_k:=\inf\{i:\ is\ge k\},\qquad k=0,1,2,\dots
$$
Let $K\in\Z$ be such that $\tau_K\le M$, $\tau_{K+1}>M$; since $Ms\ge 1$ we have $K\ge 1$. Then we have
$$
0<a_{\tau_k+1}<a_{\tau_k+2}<\dots<a_{\tau_{k+1}}\le 1\qquad\text{for each }k=0,1,2,\dots
$$
In each such increasing sequence, since $s<1/4$, we have at least $\lfloor 1/(2s)\rfloor\ge 1/(3s)$ elements in the segment $[0,1/2)$ of length $1/2$. The total number of elements in this sequence, $\tau_{k+1}-\tau_k\le \lceil 1/s\rceil\le 5/(4s)$. Consequently,
\begin{align*}
    N_{s,M}&\ge K\times \frac1{3s}
    \text{ while }
    M\le (K+1)\times \frac5{4s},
\end{align*}
so $N_{s,M}/M\ge 2/15$ since $K\ge 1$.
\end{proof}
\begin{lemma}\label{lemNEW} Let $(q_i)_{i=1}^{\infty}$ form a probability distribution ($q_i\ge 0$, $\sum_{i=1}^{\infty} q_i=1$).
Assume that for some  $a>0$ and a positive integer $M$ we have
$$
    q_j \ge \frac{a}{M}\qquad \text{for }j=1,2,\dots,M.
$$
Let $h(s):=\sum_{j=1}^\infty q_j|\cos(js)|$. Then
$$
h(s)\le 
\max\left(1-c_1' a s^2 M^2,1-c_2' a \right), \quad 0\le  s\le\frac\pi2,
$$
for some absolute constants $c_1',c_2'>0$.
\end{lemma}
\begin{proof}
First notice that 
\begin{align}\label{eq:h}
|h(s)|&\le \sum_{j=1}^M q_j|\cos(js)|+\sum_{j=M+1}^\infty q_j=
1-\sum_{j=1}^M q_j(1-|\cos(js)|)
\nonumber\\ &
\le 1-\frac{a}{M}\sum_{j=1}^M (1-|\cos(js)|).
\end{align}
We will use the elementary inequality
\begin{align}\label{eq-elemin}
    1-\cos\a\ge \frac{\a^2}{4},\quad 0\le \a\le \frac\pi2,
\end{align}
which holds since $\psi(\a):=1-\cos \a- \frac{\a^2}{4}$ has the properties $\psi(0)=0$ and $\psi'(\a)=\sin \a - \frac{\a}{2}>0$ for $\a\in(0,\pi/2)$.

\underline{Case one}: $M\ge 2$. Let $M_1=\lfloor M/2\rfloor \in[M/3, M/2]$.
If $tM_1\le \pi/2$ then by~\eqref{eq-elemin}
$$
\sum_{j=1}^M (1-|\cos(js)|)\ge \sum_{j=1}^{M_1} (1-\cos(js)) \ge \sum_{j=1}^{M_1} \frac{j^2 s^2}{4}
=\frac{M_1(M_1+1)(2M_1+1)s^2}{24}\ge \frac{M^3 s^2}{12\cdot3^3}.
$$
On the other hand, if $sM_1\ge \pi/2$ then $sM\ge 2sM_1\ge \pi$. Let $\tilde s=\frac{s}{\pi}\le \frac 12$, then  $\tilde{s} M\ge 1$ and by Lemma~\ref{lem:arithm} (setting $j_0:=-1/4$), in the set  $\{j\tilde{s} \mod 1,\ j=1,2,\dots,M\}$ there will be at least $\frac{2M}{15}$ elements which lie in the segment $[1/4,3/4)$; for the corresponding indices $j$ this implies that 
$$
|\cos(sj)|\le\max_{x\in[\frac{\pi}{4},\frac{3\pi}{4}]}
|\cos x|=\cos\frac{\pi}{4}=\frac{1}{\sqrt{2}}.
$$
Hence
$$
\sum_{j=1}^M (1-|\cos(js)|)\ge \frac{2M}{15}\left(1-\frac{1}{\sqrt 2}\right),
$$
and consequently,
$$
\frac{a}{M}\sum_{j=1}^M (1-|\cos(js)|)
\ge a\cdot \begin{cases}
\frac{M^2 s^2}{324},&\text{if }sM_1\le \pi/2;\\
\frac{2-\sqrt{2}}{15},&\text{otherwise}.
\end{cases}
$$
\underline{Case two: $M=1$.} Since $s\in [0,\pi/2]$,  \eqref{eq-elemin} yields that,
$$
\frac{a}{M}\sum_{j=1}^M (1-|\cos(js)|)
=a(1-\cos s)\ge \frac{as^2}{4}=\frac{aM^2 s^2}{4}.
$$
This, together with~\eqref{eq:h} imply the result.
\end{proof}
{\bf Acknowledgment:} We are grateful to Andrew Wade for helping us with the literature review and an anonymous referee for useful suggestions. J. E. is grateful to Lund University for its hospitality during his recent visit.

\end{document}